\documentclass[preprint,12pt]{elsarticle}
\usepackage[hyphens]{url}
\usepackage{enumerate}
\usepackage{booktabs}
\usepackage{dcolumn}
\usepackage{verbatim}
\usepackage{amsthm}
\usepackage{a4wide}
\usepackage{amsmath}
\usepackage{amssymb}
\usepackage{amstext}
\usepackage{amsthm}
\usepackage{mathrsfs}
\usepackage{graphicx}
\usepackage{psfrag}
\usepackage{url}
\usepackage{amstext}
\usepackage{multirow}
\usepackage{latexsym}
\usepackage{float}
\usepackage{subfig}
\usepackage[figuresright]{rotating} 
\usepackage{longtable}

\def\ve#1{\mathchoice{\mbox{\boldmath$\displaystyle\bf#1$}}
{\mbox{\boldmath$\textstyle\bf#1$}}
{\mbox{\boldmath$\scriptstyle\bf#1$}}
{\mbox{\boldmath$\scriptscriptstyle\bf#1$}}}

\newcommand\Z{\mathbb Z}
   
\newcommand\R{\mathbb R}
\newcommand\Q{\mathbb Q}

\DeclareMathOperator{\vol}{vol}     
\DeclareMathOperator{\Res}{Res}  

\let\epsilon=\varepsilon
\usepackage{ifthen}
\makeatletter
\newcommand{\DeclareBracket}[3]{
  \newcommand{#1}[2][]{%
  \ifthenelse%
  {\equal{##1}{}}%
  {\left#2##2\right#3}%
  {\csname ##1l\endcsname#2##2\csname ##1r\endcsname#3}}}
\makeatother
\DeclareBracket\abs||           
\DeclareBracket\norm\|\|        
\DeclareBracket\floor\lfloor\rfloor
\DeclareBracket\ceil\lceil\rceil
\DeclareBracket\set\{\}
\DeclareBracket\paren()
\DeclareBracket\inner\langle\rangle
\DeclareBracket\fractional\{\}

\newcommand\C{\mathbb C}

\makeatletter
\newtheorem{lemma}{Lemma}
\renewcommand*{\c@lemma}{\c@theorem}
\renewcommand*{\p@lemma}{\p@theorem}

\renewcommand*{\c@conjecture}{\c@theorem}
\renewcommand*{\p@conjecture}{\p@theorem}

\newtheorem{proposition}{Proposition}
\renewcommand*{\c@proposition}{\c@theorem}
\renewcommand*{\p@proposition}{\p@theorem}

\newtheorem{corollary}{Corollary}
\renewcommand*{\c@corollary}{\c@theorem}
\renewcommand*{\p@corollary}{\p@theorem}

\renewcommand*{\c@observation}{\c@theorem}
\renewcommand*{\p@observation}{\p@theorem}

\theoremstyle{definition}

\renewcommand*{\c@problem}{\c@theorem}
\renewcommand*{\p@problem}{\p@theorem}

\renewcommand*{\c@definition}{\c@theorem}
\renewcommand*{\p@definition}{\p@theorem}

\renewcommand*{\c@remark}{\c@theorem}
\renewcommand*{\p@remark}{\p@theorem}

\renewcommand*{\c@example}{\c@theorem}
\renewcommand*{\p@example}{\p@theorem}

\makeatother

\usepackage[breaklinks=true,colorlinks,citecolor=blue]{hyperref}

\usepackage{algorithm}
\usepackage{algorithmic}  
\renewcommand{\algorithmicrequire}{\textbf{Input:}}
\renewcommand{\algorithmicensure}{\textbf{Output:}}

\newcommand{\timeFaster}[1]{{\bf #1}}
\newcommand{\timeSlower}[1]{#1}
\newcommand{\timeTie}[1]{#1}
\newcommand{\timeNotComputed}[1]{--}

\usepackage{xspace}

\newcommand{\latte}{{\tt LattE}\xspace}
\newcommand{\latteInt}{{\tt LattE integrale}\xspace}

\newcommand{\vinci}{{\tt Vinci}\xspace}
\newcommand{\maple}{{\tt Maple}\xspace}
\newcommand{\cdd}{{\tt Cddlib}\xspace}
\newcommand{\fourtitwo}{{\tt 4ti2}\xspace}
\newcommand{\cubpack}{{\tt CUBPACK}\xspace}

\renewcommand\d{\,\mathrm{d}}

\DeclareBracket\charfun[]       
\newcommand{\e}{{\mathrm{e}}}   

\journal{Computational Geometry: Theory and Applications}

\begin{document}

\begin{frontmatter}
\title{Software for Exact Integration of Polynomials over Polyhedra}

\title{Software for Exact Integration of Polynomials over Polyhedra}
\author[add1]{J. A. De Loera}
\ead{deloera@math.ucdavis.edu}

\author[add1]{B. Dutra\corref{cor1}}
\ead{bedutra@ucdavis.edu}

\author[add1]{M. K\"oppe}
\ead{mkoeppe@math.ucdavis.edu}

\author[add1]{S. Moreinis\fnref{stanford}}
\ead{moreinis@stanford.edu}

\author[add1]{G. Pinto}
\ead{gpinto@ucdavis.edu}

\author[add1]{J. Wu\fnref{caltech}}
\ead{jwu2@caltech.edu}

\address[add1]{Department of Mathematics, UC Davis, One Shields Avenue, Davis, CA 95616, USA}

\fntext[stanford]{Present address: Stanford University, Stanford, California 94305}
\fntext[caltech]{Present address: California Institute of Technology, Pasadena, California 91125}
 \cortext[cor1]{Corresponding author.}

\begin{abstract} 
We are interested in the fast computation of the exact value of
integrals of polynomial functions over convex polyhedra. We
present speed ups and extensions of the algorithms
presented in previous work by some of the authors. We present a new software
implementation and provide benchmark computations.  The
computation of integrals of polynomials over polyhedral regions has
many applications; here we demonstrate our algorithmic tools solving
a challenge from combinatorial voting theory.
\end{abstract}


\end{frontmatter}


{\footnotesize\tableofcontents}

\section{Introduction}

Integration is a fundamental problem arising in many contexts, from
engineering to algebraic geometry. For example,  center of mass, moments and products of
inertia about various axes are integrals used in any type of dynamic
simulation or modeling such as computer games
\cite{BLi:1993,Mirtich:1996, Ong-2003:finding-exact-volume-polyhedron}), similarly,
normalized volumes indicate the degrees of toric varieties and closely related to moment maps of symplectic
manifolds \cite{Sturmfels1996,Berline-Vergne-2007}). Integration over
polyhedra is particularly useful because many domains can be
approximated by polyhedra and then decomposed into convex polyhedra
(e.g., a tetrahedral mesh decomposition etc.).

In this work we are interested in the exact \emph{evaluation
  of integrals} over convex polyhedral regions.  More precisely, let
$P$ be a $d$-dimensional rational convex polyhedron in $\R^n$ and
let $f\in \Q[x_1,\dots,x_n]$ be a polynomial with rational
coefficients. We consider the problem of efficiently computing the
\emph{exact} value of the integral of the polynomial $f$ over $P$,
which we denote by $\int_{P} f \d m$. Here we use the \emph{integral
  Lebesgue measure} $\d m$ on the affine hull $\inner{ P}$ of the
polytope $P$. This general setting is quite useful because, when the
polytope is full-dimensional, the integral Lebesgue measure coincides
with the standard Riemann integral but generalizes it naturally to
cases when the polytope is not full-dimensional. Another reason to
compute in this setting is that the volume of $P$ and every integral
of a polynomial function with rational coefficients yields
\emph{rational numbers}.  Finally this normalization of the measure
occurs naturally in Euler--Maclaurin formulas for a polytope $P$,
which relate sums over the lattice points of~$P$ with integrals over
the various faces of~$P$. We remark that 
the computer algebra community has dedicated a great deal of effort to
develop a different kind of exact integration, understood to be the automatic
computation of the antiderivatives of functions, as predicted by the
fundamental theorem of Calculus \cite{Bronstein}, but  we do not discuss this
kind of exact integration here.

Regarding the theoretical computational complexity of our problem, 
it is very educational to look first at the case when $f$ is the constant
polynomial~$1$, and the answer is simply a volume.  It has been
proved that already computing the volume of polytopes of varying
dimension is $\#\mathrm{P}$-hard 
\cite{dyerfrieze88,brightwellwinkler91,khachiyan93,lawrence91, gritzmann-complexity-volume:1994}, and
that even approximating the volume is hard~\cite{elekes86}. More 
recently in \cite{rademacher} it was proved that computing the
center of mass of a polytope is $\#\mathrm{P}$-hard.

We report on a new C++ implementation, {\tt LattE integrale} \cite{latteintegrale}, which extends the 
work done in \cite{baldoni-berline-deloera-koeppe-vergne:integration} and \cite{bueler-enge-fukuda-2000:exact-volume}.
This paper is mostly an experimental and practical study, but it also slightly  develops the theory of 
\cite{baldoni-berline-deloera-koeppe-vergne:integration}. This article presents useful
formulas for integration of power of linear forms over simplicial cones that complement
those presented in \cite{baldoni-berline-deloera-koeppe-vergne:integration}.

Our method of computation relies on powerful mathematical ideas.  It
was  proved in \cite{baldoni-berline-deloera-koeppe-vergne:integration}
that, unlike general polynomials, integrals over simplices of
arbitrary powers of linear forms, or of polynomials with a fixed number
of variables, can be computed in polynomial time. In this case our
algorithms use known properties of integrals of exponentials of linear forms
(see \cite{Barvinok-1991}, \cite{Barvinok-1992}). This allows for fast calculation 
over general polytopes using two methods that depend on two different decompositions
of polyhedra. General polyhedra can be decomposed as either a disjoint union of simplices, 
i.e., triangulations, or as signed cone decompositions of the kind proposed by Brion, Lasserre, 
Lawrence, and Varchenko \cite{Brion88, lasserre-algorithm-volume-polyhedron:1983,lawrence91, brion-lawrence-varchenko-rational-generating-functions-2005}. The polynomial-time complexity for integration over simplices shown in \cite{ baldoni-berline-deloera-koeppe-vergne:integration} can be extended to more polyhedra as long as their decompositions 
are of ``small'' size  (note that this is always the case in fixed dimension).

This paper is organized as follows. We begin in Section~\ref{preliminar} recalling 
the mathematical ideas at the heart of our algorithms (although we omit details of proofs, they can be
found in the references). We begin with a short review of polyhedral geometry, specially valuations. In
Section~\ref{section-software} we discuss details about the
implementation including main subroutines and data structures.  In
Section~\ref{section-integration-experiments} we first discuss speed
improvements for integrating over simplices from earlier work in
\cite{baldoni-berline-deloera-koeppe-vergne:integration}, and then we
report on several benchmarks of integration over arbitrary polytopes. More experimental tables
are available online \cite{latteintegrale}. We conclude our paper with an application: we solve a 
computational challenge from combinatorial voting theory.

\section{Mathematical preliminaries} \label{preliminar}

In this section we recall the necessary mathematical background used in our algorithms. 
We state results without proof but excellent background sources for what is going to be discussed here include
 \cite{baldoni-berline-deloera-koeppe-vergne:integration,barvinokzurichbook,BarviPom,beckrobins} and the
references mentioned there.

\subsection{Polyhedra and polynomials}
A \emph{convex rational polyhedron} $P$ in $\R^d$ (we will simply say
 \emph{polyhedron}) is defined as the intersection of a finite number $m$ of
closed half spaces bounded by rational affine hyperplanes. We say
that  $P$ is \emph{full-dimensional} (in $\R^d$) if the affine span of $P$ is $\R^d$.
When $P=\{\, x : Ax \leq b \, \}$ for a $m \times d$ matrix $A$ and $m$-vector $b$, $P$ 
is said to be given by a \emph{halfspace} or \emph{h-representation}. When $P$ is 
the convex hull conv($V$) of finitely many points in $\R^d$, $V = \{v_1, \dots v_n\}$, 
$P$ is said to be given by a \emph{vertex} or \emph{v-representation}. We can switch 
between the h- and v-representations of a $d$-dimensional polyhedron using
well-known algorithms (see \cite{Avis-lrs-1999, fukuda-cdd-revisited:1996}).
A \emph{polytope} $P$ is a  compact polyhedron.
A \emph{cone} $C$ is a polyhedral cone (with vertex $0$) and
an \emph{affine cone} is a translation $s+C$ of a cone $C$. A  cone
$C$ is called \emph{simplicial} if it is  generated by linearly independent
vectors of $\R ^d$. A simplicial cone~$C$ is called \emph{unimodular} if it
is generated by  linearly independent integral vectors $v_1,\dots, v_k$ such
that $\{v_1,\dots, v_k\}$ can be completed to an integral basis of
$\Z^d$. An affine cone $C$ is called simplicial (respectively, simplicial
unimodular) if the associated cone is.  
 The set of vertices of $P$ is denoted by $V(P)$.
For each vertex $s\in V(P)$, the \emph{cone of feasible directions} $C_s (P)$ 
at the vertex $s$  is the cone of all vector $y$ such that $s+\epsilon y \in P$ for some $\epsilon>0$.
 The \emph{tangent cone} of a polytope $P$ at a vertex $s$ 
is the affine cone $s+C_s(P)$ (this is a translation of $C_s(P)$). For details in all 
these notions see, e.g., \cite{barvinokzurichbook}.

For the integration of a full-dimensional polytope we consider the standard 
Lebesgue measure, which gives volume $1$ to the fundamental domain
of the lattice $\Z^n$. But if $P$ lies inside an affine subspace $L+a$, with $L$ a rational linear subspace of
dimension $n\leq d$, we will normalize the Lebesgue measure on~$L$, so
that the volume of the fundamental domain of the lattice
$L\cap \Z^d$ is~$1$. Thus for any affine subspace $L+  a$
parallel to $L$, we define the \emph{integral Lebesgue measure}  $\d
m$  by translation. For example, the diagonal  of the unit square
has length $1$ instead of  $\sqrt{2}$. This has the great advantage that for
rational input our output will always be an (exact) rational number
$\int_P f \d m$ in the usual binary encoding. 

 One important point of our method is that all computations are done 
 in the representation polynomials given by \emph{powers of linear forms}. It is 
 well-known that any homogeneous polynomial of degree $M$ can be
 decomposed as a sum of $M$-th powers of linear forms. For example, 
one can decompose the polynomial $f$ as a sum 
$f =\sum_{\ell}c_{\ell} \langle \ell, x \rangle^M$ with at most $2^M$ terms. This decomposition
is given by the following well-known identity for monomials:
If $x^{\ve M}= x_1^{M_1}x_2^{M_2}\cdots x_n^{M_n}$, then

\begin{multline} 
x^{\ve M}  = \frac{1}{|\ve M|!} \sum_{0\leq p_i\leq M_i}(-1)^{|\ve M|-(p_1+\cdots+p_n)}
  \binom{M_1}{p_1}\cdots \binom{M_n}{p_n}(p_1 x_1+\cdots+p_n x_n)^{|\ve M|},
  \label{eq:decomp-powerlinform}
\end{multline}
where $|\ve M| = M_1+\cdots+M_n \leq M$. Of course, when dealing with 
general polynomials, this same formula can be applied for as many monomials as is necessary.
For example,  the polynomial $7x^2 + y^2 + 5z^2 + 2xy + 9yz$ can be written as
 $\frac{1}{8} (12(2x)^2 - 9(2y)^2 + (2z)^2 + 8(x + y)^2 + 36(y + z)^2).$

It is worth noting that the above formula does not yield an optimal decomposition, 
but it suffices to generate a polynomial-time algorithm on fixed degree $|\ve M|$ or 
fixed number of variables~\cite{baldoni-berline-deloera-koeppe-vergne:integration}. 
The problem of finding a decomposition with the smallest possible number of summands 
is known as the \emph{polynomial Waring problem}.  What is the smallest integer $r(M, n)$ 
such that a generic homogeneous polynomial $f(x_1, \dots ,x_n)$ of degree~$M$ in $n$ 
variables is expressible as the sum of $r(M,n)$ $M$-th powers of linear forms?  
This problem was solved for general polynomials by Alexander and Hirschowitz \cite{alexanderhirschowitz} 
(see \cite{brambillaottaviani} for an extensive survey), but there is no computational or constructive version of 
this result that would yield the optimal decomposition for an specific input polynomial and
the  bounds may be much too pessimistic on the average situation. Only very recently Carlini et al. \cite{Carlinietal} 
gave efficient decompositions of a monomial.  However, their decomposition involves roots of unity, 
and here we are interested in an arithmetic version of the problem where everything is expressed using 
rational forms and rational coefficients. But we can see that the explicit rational construction we use in our
code is not too far away from the optimum.

Table~\ref{tabel:count-linear-forms} lists the average number of powers of linear forms necessary to 
decompose monomials of given degree generated uniformly at random. To create the monomials, we keep 
adding one to the power of a randomly chosen variable until the monomial has the desired degree. The
table show mild exponential growth as degree or dimension grow. This was predicted in the theory.

\begin{table}
\centering
\caption{Average number of powers of linear forms plus or minus one standard deviation necessary to express one monomial in $d$ variables, averaged over 50 monomials of the same degree}
\label{tabel:count-linear-forms}
\tabcolsep 2.5pt
\small
\begin{tabular}{c*8{r}}
	\toprule	
	
	& \multicolumn{6}{c}{Monomial Degree}\\
		\cmidrule(r){2-7}

\multicolumn{1}{c}{$d$}&   \multicolumn{1}{c}{5} &  \multicolumn{1}{c}{10} &  \multicolumn{1}{c}{20} &  \multicolumn{1}{c}{30} &  \multicolumn{1}{c}{40} &  \multicolumn{1}{c}{50} \\ 
\hline
3 & $14 \pm 3$ & $(6.6 \pm 1.2)\times10^{1}$ & $(4.0 \pm 0.5)\times10^{2}$ & $(1.2 \pm 0.1)\times10^{3}$ & $(2.7 \pm 0.2)\times10^{3}$ & $(5.2 \pm 0.2)\times10^{3}$  \\ 
4 & $16 \pm 5$ & $(1.1 \pm 0.2)\times10^{2}$ & $(1.1 \pm 0.2)\times10^{3}$ & $(4.5 \pm 0.6)\times10^{3}$ & $(1.3 \pm 0.2)\times10^{4}$ & $(3.0 \pm 0.2)\times10^{4}$  \\ 
5 & $19 \pm 4$ & $(1.5 \pm 0.4)\times10^{2}$ & $(2.2 \pm 0.6)\times10^{3}$ & $(1.2 \pm 0.3)\times10^{4}$ & $(4.7 \pm 0.7)\times10^{4}$ & $(1.4 \pm 0.2)\times10^{5}$  \\ 
6 & $20 \pm 5$ & $(2.0 \pm 0.6)\times10^{2}$ & $(4.1 \pm 1.2)\times10^{3}$ & $(3.2 \pm 0.8)\times10^{4}$ & $(1.5 \pm 0.3)\times10^{5}$ & $(5.2 \pm 0.6)\times10^{5}$  \\ 
7 & $21 \pm 5$ & $(2.4 \pm 0.9)\times10^{2}$ & $(6.7 \pm 2.4)\times10^{3}$ & $(7.1 \pm 2.1)\times10^{4}$ & $(4.0 \pm 1.0)\times10^{5}$ & $(1.7 \pm 0.3)\times10^{6}$  \\ 
8 & $21 \pm 5$ & $(2.9 \pm 0.9)\times10^{2}$ & $(1.1 \pm 0.5)\times10^{4}$ & $(1.4 \pm 0.5)\times10^{5}$ & $(9.8 \pm 2.7)\times10^{5}$ & $(4.8 \pm 1.1)\times10^{6}$  \\ 
10 & $24 \pm 5$ & $(3.5 \pm 1.1)\times10^{2}$ & $(2.1 \pm 0.9)\times10^{4}$ & $(4.1 \pm 1.6)\times10^{5}$ & $(4.5 \pm 1.7)\times10^{6}$ & $(3.1 \pm 1.0)\times10^{7}$ \\ 

	\bottomrule
\end{tabular}
\end{table}


In conclusion, to integrate a multivariate polynomial, we first algebraically decompose each monomial 
to a sum of powers of linear forms which, as we will see next, can be integrated very fast in practice
over simplices or over simplicial cones using a few useful formulas. Thus we will need a geometric 
decomposition of our polytopes into those pieces.

\subsection{Valuations and formulas of integration of exponentials over cones and simplices}

We now recall  several formulas  for the integrals of a power of a linear form over
a simplex or over a simplicial cone. The idea is that if we can do  fast integration
for those two structures, then we can always rely on two polyhedral decompositions of the input 
polyhedron to obtain the integral. See Subsection~\ref{section-triangulate-or-decompose} for details.
 
One of the most important properties of integrals over polyhedra is that
they can be seen as valuations.  A \emph{valuation} $F$ is a linear map 
from the rational vector space of the indicator functions of rational
polyhedra $ P\subset \R^d$ into a rational vector space $M$. Whenever
the indicator functions $\charfun {P_i}$ of a family of polyhedra
$P_i$ satisfy a linear relation $\sum_i r_i \charfun{P_i}=0$, then the
elements $F(P_i)$ satisfy the same relation $ \sum_i r_i F(P_i)=0$
(for a formal definition within the polytope algebra, see Chapter 2 of  \cite{barvinokzurichbook}).

Let $C=\sum_{i=1}^d \R_+ u_i$ be the simplicial cone spanned by linearly independent integral vectors $u_1,u_2,\dots u_d$.
The \emph{fundamental parallelepiped} $\Pi_C$ of the cone $C$
(with respect to the generators $u_i,i=1, \dots,d$) is the set of
points $\Pi_C=\sum_{i=1}^d [0,1\mathclose[\, u_i$. Let us denote by
    $\vol(\Pi_C)$ its volume.

\begin{proposition}[Theorem 8.4 in \cite{barvinokzurichbook}].\label{valuationI}
There exists a unique valuation  $I(P)(\ell)$ which  associates  to every polyhedron
$P\subset V$ a meromorphic function so that the following properties hold 

\begin{enumerate}[\rm(i)]

\item If $\ell$ is a linear form such that $\e^{\langle \ell, x\rangle}$ is integrable over $P$,
  then
$$
I(P)(\ell)= \int_P \e^{\langle \ell,x\rangle} \,\mathrm dm.
$$

\item For every point $s\in \R^n$, one has
$$
I(s+P)(\ell) = \e^{\langle \ell,s\rangle}I(P)(\ell).
$$
\item If $P$ contains a straight line, then $I(P)=0$.
\end{enumerate}
\end{proposition}

A consequence of the valuation property is the following fundamental theorem. 
It follows from the Brion--Lasserre--Lawrence--Varchenko decomposition theory of a
polyhedron into the supporting cones at its vertices \cite{Brion88,barvinokzurichbook, brion-lawrence-varchenko-rational-generating-functions-2005, lasserre-algorithm-volume-polyhedron:1983}.

\begin{lemma} \label{brion-exp} Let $P$ be a polyhedron with set of vertices $V(P)$. For each
vertex~$s$, let $C_s(P)$ be the cone of feasible directions at vertex $s$. Then
\begin{equation}
I(P)(\ell)=\sum_{s\in V(P)}I(s+C_s(P))(\ell).
\end{equation}
\end{lemma}

Note that the cone $C_s(P)$ in Lemma~\ref{brion-exp} may not be simplicial, but for simplicial cones their integrals 
have explicit rational function formulas. As we see in Proposition~\ref{prop:integral-exp-simplicial}, one can derive an explicit formula for 
the rational function $I(s+C_s(P))$ in terms of the geometry of the cones.

\begin{lemma} \label{pieces}
Using the valuation property for the valuation $I(P)(\ell)$ and the linearity over the integrands we have that:

\begin{enumerate}[\rm(i)]
\item For any triangulation $\mathcal T$ of the polytope $P$, we have
  $I(P)(\ell)=\sum_{\Delta \in \mathcal T} I(\Delta)(\ell)$.

\item For any triangulation $\mathcal D_s$ of the feasible cone~$C_s(P)$ at each of the
  vertices~$s$ of the polytope~$P$ we have $I(P)(\ell)=\sum_{s\in V(P)} \sum_{C
    \in \mathcal D_s} I(s+C)(\ell).$ 
\end{enumerate}
\end{lemma}

Lemma~\ref{pieces} says that \emph{if} we know how to integrate over simplices or 
simplicial cones, we can integrate over a polytope. We are close to knowing how to do this. By
elementary integration, and Proposition~\ref{valuationI}, we have the following. 

\begin{proposition} 
  \label{prop:integral-exp-simplicial}
  For a simplicial cone $C$ generated by rays $u_1,u_2,\dots u_d$ (with vertex $0$) and for any point $s$
\begin{equation}
I(s+C)(\ell)
=\vol(\Pi_C) \e^{\langle \ell, s
  \rangle} \prod_{i=1}^d \frac1{\langle-\ell, u_i \rangle}.
\end{equation}
This identity holds as a meromorphic function of~$\ell$ 
and pointwise for every $\ell$ such that $\langle \ell, u_i \rangle \neq 0$ for
all $u_i$.
\end{proposition}

\subsection{From exponentials to powers of linear forms}

We now consider powers of linear forms instead of exponentials.  Similar to
$I(P)$, we now let $L^M(P)$ be the meromorphic extension of the function defined by
$$ L^M(P)(\ell) = \int_{P} {\langle \ell, x \rangle}^M  \d{m}$$ for those
$\ell$ such that the integral exists. 
To transfer what we know about integrals of exponentials to those of powers of
linear forms, 
we can consider the formula of Proposition~\ref{prop:integral-exp-simplicial}
as a function of the auxiliary parameter $t$:  
\begin{equation}\label{eq:formula-with-t}
  \int_{s+C} \e^{\langle t \ell, x \rangle}  \d{m}=\vol(\Pi_C)
  \e^{\langle t \ell, s \rangle}\prod_{i=1}^d \frac{1}{\langle- t \ell,
    u_i \rangle}. 
\end{equation}
Using the series expansion of the left in the variable~$t$, we wish to
recover the value of the integral of $\langle \ell, x \rangle^M$ over
the cone. This is the coefficient of $t^M$ in the expansion; to
compute it, we equate it to the Laurent series expansion around $t=0$
of the right-hand-side expression, which is a meromorphic function
of~$t$. Clearly
$$\vol(\Pi_C) \e^{\langle t \ell, s \rangle}  \prod_{i=1}^d \frac1{\langle- t
  \ell, u_i \rangle}=\sum_{n=0}^\infty t^{n-d} \frac{\langle \ell, s
  \rangle^n}{n!}\cdot \vol(\Pi_C) \prod_{i=1}^d\frac1{\langle-\ell, u_i \rangle}.$$ 
We say that $\ell$ is \emph{regular}  if $\langle \ell, u_i \rangle \neq 0$ for 
 every ray $u_i$ of the cone. With this, we can conclude the following.

\begin{corollary} \label{oneconeoneplf}
  For a regular linear form $\ell$ and a simplicial cone $C$ generated by rays $u_1,u_2,\dots u_d$ with vertex $s$
\begin{equation} \label{above}
L^M(s+C)(\ell) 
=\frac{M!}{(M+d)!} \vol(\Pi_C) \frac{(\langle \ell, s \rangle)^{M+d}}{\prod_{i=1}^d  \langle -\ell, u_i \rangle}.
\end{equation}
Otherwise when $\ell$ is not regular,  there is a nearby perturbation which is regular. To obtain it,
we use  $\ell + \hat \epsilon$ where $\hat \epsilon = \epsilon a$ is any linear form with $a \in \R^n$ 
such that $\langle -\ell - \hat \epsilon, u_i \rangle \neq 0$ for all $u_i$,  to define a new linear form 
(depending of $a$) on the space of meromorphic functions in the variable $\epsilon$. Then, applying \eqref{above} on 
the limit as $\epsilon$ goes to zero we obtain:
\begin{equation}\label{eq:residue-with-form-perturbation}
L^M(s+C)(\ell) 
=\frac{M!}{(M+d)!} \vol(\Pi_C) \Res_{\epsilon=0}\frac{(\langle \ell + \hat \epsilon, s \rangle)^{M+d}}{\epsilon \prod_{i=1}^d  \langle - \ell - \hat \epsilon, u_i \rangle},
\end{equation}
\end{corollary}

Since the reader may not be familiar  with residues $\Res_{\epsilon=0} (f)$ and how to calculate them,  we recall some useful facts
on complex analysis (see, e.g., \cite{Henrici} for details).  As we observed, there is a  singularity or pole at $\epsilon=0$ for a univariate rational 
function $f(\epsilon)=\frac{p(\epsilon)}{q(\epsilon)}$ (which in this case is explicitly given in Formula \eqref{eq:residue-with-form-perturbation} of Corollary \ref{oneconeoneplf}). 
Recall that if $f$ has Laurent expansion $f(\epsilon)=\sum_{k=-m}^{\infty} a_k \epsilon^k,$ the residue is defined as $a_{-1}$.
Given a rational function $f(\epsilon)$ with a pole at~$\epsilon=0$ there are a variety of well-known techniques to extract the value of the residue. 
For example,    if $\epsilon=0$ is a simple pole ($m=1$), then
$\Res_{\epsilon=0} (f)=\frac{p(0)}{q'(0)}$. 
Otherwise, when $\epsilon=0$ is a pole of order $m>1$, we can write $f(\epsilon)=\frac{p(\epsilon)}{\epsilon^m q_1(\epsilon)}.$
then  expand $p,q_1$ in powers of $\epsilon$  with $p(\epsilon)=a_0+a_1\epsilon+a_2\epsilon^2+\dots$ and $q_1(\epsilon)=b_0+b_1\epsilon+b_2\epsilon^2+\dots$.
This way the Taylor expansion of $p(\epsilon)/q_1(\epsilon)$ at $\epsilon_0$ is $c_0+c_1\epsilon+c_2\epsilon^2+c_3\epsilon^3+\dots$, where
$ c_0=\frac{a_0}{b_0}$, and
$c_k=\frac{1}{b_0}(a_k-b_1c_{k-1}-b_2c_{k-2}-\dots-b_kc_0)$. Thus  we recover
the residue  
$\Res_{\epsilon=0}(f)=c_{m-1}$. We must stress that the special structure of the rational functions in Corollary \ref{oneconeoneplf} can be exploited to
speed up computation further rather than using this general methodology. For more on this see \cite{barvinok-woods-2003, barvinokzurichbook, baldoni-berline-deloera-koeppe-vergne:integration} 
and the following discussion.

Finally, we have all the tools necessary to write down our formula for integration using cone decompositions.

\begin{corollary} \label{cor:the-lawrence-method}
For any triangulation $\mathcal D_s$ of the tangent cone~$C_s(P)$ at each of the
vertices~$s$ of the polytope~$P$ we have 
\begin{equation}
	L^M(P)(\ell)
        =\sum_{s\in V(P)} \sum_{C \in \mathcal D_s} L^M(s+C)(\ell).
\end{equation}
\end{corollary}

\subsection{The formula for the simplex}

Suppose now  that $\Delta \subset \R^n$ is a $d$-dimensional simplex (as it
may appear in a triangulation of the polytope~$P$), and
$\ell$ is a linear form on $\R^n$.  We say that the linear form~$\ell$ is
\emph{regular} for the simplex~$\Delta$ if it is not orthogonal to any of the edges of
the simplex.  If $\ell$ is regular for~$\Delta$, then it is regular for all
tangent cones at each of the vertices.  We then find the following result
as a special case of Corollary~\ref{cor:the-lawrence-method}.

\begin{corollary}[Brion, see \cite{Brion88}] \label{th:brion}
  Let $\Delta$ be a d-simplex with vertices $s_1, \dots, s_{d+1} \in \R^d$.  Let $\ell$
  be a linear form which is regular w.r.t.~$\Delta$, i.e., $\langle
  \ell, s_i\rangle \neq \langle \ell, s_j \rangle$ for any
  pair $i\neq j$. Then we have the following relation.
\begin{equation}\label{eq:brion-powerlinform-regular}
  L^M(\Delta)(\ell)
  =\int_\Delta {\langle \ell, x \rangle}^M \mathrm d m = d!\vol(\Delta, \d m)\frac{M!}{(M+d)!}
  \Big(\sum_{i=1}^{d+1}\frac{ \langle \ell , s_i
    \rangle^{M+d}}{\prod_{j\neq i} \langle \ell, s_i- s_j \rangle}\Big).
\end{equation}
\end{corollary}
When $\ell$ is regular, Brion's formula is very short; it is a sum
of $d+1$ terms. When $\ell$ is not regular, we can again use a perturbation
$\ell + \hat \epsilon$ where $\hat \epsilon = \epsilon a$ as in
Corollary~\ref{oneconeoneplf}, so that the expression of  the integral over
the simplex reduces to a sum of residues as
in~\eqref{eq:residue-with-form-perturbation}. 

However, in the special case of a simplex, there exists a computationally more
efficient method that avoids the calculation of a perturbation~$a$; see
\cite{baldoni-berline-deloera-koeppe-vergne:integration}.  From
\cite[Theorem~10]{baldoni-berline-deloera-koeppe-vergne:integration} we find
that $L^M(\Delta)(\ell)$ is the coefficient of the term~$t^M$ in the Laurent
series of the rational function
\begin{equation}
  \label{eq:howto-paper-theorem10}
  d! \vol(\Delta, \d m)\frac{M!}{(M+d)!} \frac{1}{\prod_{j=1}^{d+1}(1-t\langle \ell , s_j\rangle)}
\end{equation}
in the variable~$t\in\C$.  This rational function can be expanded into partial
fractions. To this end,
let $K\subseteq\{1,\dots,d+1\}$ be an index set of the different poles
$t= t_k := 1/\langle \ell ,s_k\rangle$, and for $k\in K$ let $m_k$ denote
the order of the pole, i.e.,
\begin{displaymath}
  m_k = \#\bigl\{\, i\in\{1,\dots,d+1\} : \langle \ell ,s_i\rangle = \langle \ell ,s_k\rangle \,\bigr\}.
\end{displaymath}
Then the rational function can be written as 
\begin{displaymath}
  \sum_{k\in K} \Bigl(\frac{a_{k,1}}{1-t\langle \ell , s_k\rangle} +
  \frac{a_{k,2}}{(1-t\langle \ell , s_k\rangle)^2} + \dots +
  \frac{a_{k,m_k}}{(1-t\langle \ell , s_k\rangle)^{m_k}} \Bigr),
\end{displaymath}
where the coefficients $a_{k,r}$ are given by certain residues about the pole $t=t_k$. 
After a change of variables, $t = t_k + \epsilon$, one obtains the following formula.

\begin{corollary}[Corollary 13 in~\cite{{baldoni-berline-deloera-koeppe-vergne:integration}}]\label{th:residue}
Let $\Delta$ be a $d$-dimensional simplex. Then for an arbitrary power
$\langle \ell, x \rangle^M$ of a linear form, we have:
\begin{equation}
  \int_{\Delta} {\langle \ell, x \rangle}^M  \d m =
 d!\vol(\Delta, \d m) \frac{M!}{(M+d)!}\sum_{k\in K} \Res_{\epsilon=0} \frac{(\epsilon + \langle
    \ell,  s_k \rangle)^{M+d}}
  {\epsilon^{m_k} {\prod\limits_{\substack{i\in K\\ i\neq k}} {(\epsilon +
      \langle \ell,  s_k-s_i\rangle )}^{m_i}} }.
  \label{eq:integral-via-residues}
\end{equation}

\end{corollary}

To conclude we note that one can even extend the formula above on
integrating the power of a linear form to the case of a \emph{product
  of powers of several linear forms} (see
\cite{baldoni-berline-deloera-koeppe-vergne:integration}).

\subsection{Should one triangulate or cone decompose?}
\label{section-triangulate-or-decompose}

One could triangulate the whole polytope, or integrate over each
tangent cone.  However, each cone must be decomposed into simplicial
cones.  This is the trade-off: we can get away with not doing one
large polytope triangulation, but we might have to do many smaller
cone triangulations.

The number of simplices in a triangulation and the number of
simplicial cones containing in a polytope can significantly differ. Depending on
the polytope, choosing the right method can determine its
practicality. Our experimental results agree with
\cite{bueler-enge-fukuda-2000:exact-volume} in showing that triangulating the
polytope is better for polytopes that are ``almost simplicial'' while cone
decomposition is faster for simple polytopes. The details will be discussed
in Section~\ref{section-integration-experiments}.

\bigbreak
\noindent Lemma~\ref{pieces} together with the formulas we stated for integration over simplices and cones give a general process for computing integrals:
\begin{itemize} 
\item We decompose our polynomial as a sum of powers of linear forms.
\item We select a decomposition of the polyhedron in question, either a triangulation or a cone decomposition.
\item Apply the formulas to each piece and add up the results via the above results.
\end{itemize}
\subsection{Examples}
\label{section-examples-of-each}

\subsubsection{Integral values encoded by rational function identities}
\begin{figure}
 	\centering
	\includegraphics[width=0.3\textwidth]{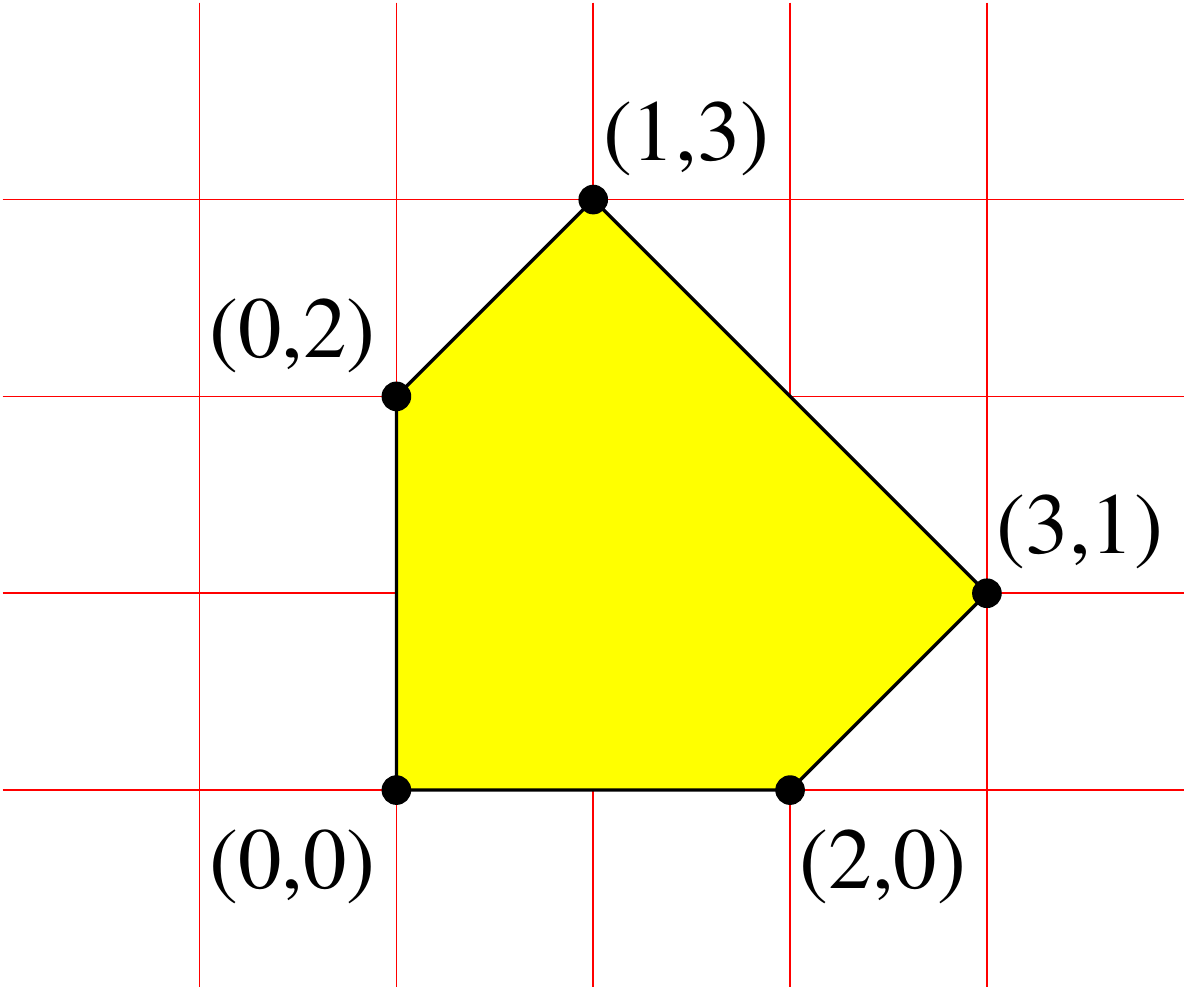}                
	\caption{A pentagon}
	\label{fig:pentagon}
\end{figure}

Before working out a simple integration example, let us highlight the
fact that for regular linear forms the integration formulas are given
by sums of rational functions which we read from the geometry at vertices
and possibly a cone decomposition method: Consider a pentagon $P$ with
vertices $(0,0), (2,0), (0,2), (3,1),$ and $(1, 3)$ as in Figure
\ref{fig:pentagon}.

Then the rational function giving the value of $\int_{P} (c_1 x +c_2 y)^{M} \d{x} \d{y}$ is
\begin{displaymath}
	\frac{M!}{(M+2)!} \, \left( {\frac { \left( 2\,c_{{1}} \right) ^{M+2}}{c_{{1}} \left( 
-c_{{1}}-c_{{2}} \right) }}+4\,{\frac { \left( 3\,c_{{1}}+c_{{2}}
 \right) ^{M+2}}{ \left( c_{{1}}+c_{{2}} \right)  \left( 2\,c_{{1}}-2
\,c_{{2}} \right) }}+4\,{\frac { \left( c_{{1}}+3\,c_{{2}} \right) ^{M
+2}}{ \left( c_{{1}}+c_{{2}} \right)  \left( -2\,c_{{1}}+2\,c_{{2}}
 \right) }}+{\frac { \left( 2\,c_{{2}} \right) ^{M+2}}{ \left( -c_{{1}
}-c_{{2}} \right) c_{{2}}}} \right).
\end{displaymath}

This rational function expression encodes \emph{every integral} of the form $\int_{P} (c_1 x +c_2 y)^{M} \d{x} \d{y}$. For example, if we let $M=0$, then the integral is equal to the area of the pentagon, and the rational function simplifies to a number by simple high-school algebra:
\begin{displaymath}
\frac{1}{2} \left( 4\,{\frac {c_{{1}}}{-c_{{1}}-c_{{2}}}}+4\,{\frac { \left( 3\,c_{{1}}+c
_{{2}} \right) ^{2}}{ \left( c_{{1}}+c_{{2}} \right)  \left( 2\,c_{{1}
}-2\,c_{{2}} \right) }}+4\,{\frac { \left( c_{{1}}+3\,c_{{2}} \right) 
^{2}}{ \left( c_{{1}}+c_{{2}} \right)  \left( -2\,c_{{1}}+2\,c_{{2}}
 \right) }}+4\,{\frac {c_{{2}}}{-c_{{1}}-c_{{2}}}} \right)  = 6.
 \end{displaymath}

Hence the area is $6$. When $M$ and $(c_1,c_2)$ are given and $(c_1,c_2)$ is not perpendicular to any of the edge directions we can simply plug in numbers to the rational function. 
For instance, when $M=100$ and $(c_1=3,c_2=5)$ the answer is a fraction with numerator equal to
\begin{align*}
227276369386899663893588867403220233833167842959382265474194585 \\
3115019517044815807828554973991981183769557979672803164125396992
\end{align*}
and denominator equal to $1717$. When $(c_1,c_2)$ is perpendicular to
an edge direction, we encounter (removable) singularities in the rational functions,
thus using complex residues we can do the evaluation.  Note that those linear forms that 
are perpendicular to some edge direction form a measure zero set inside
a hyperplane arrangement.

\subsubsection{Using the triangulation method}
\begin{figure}
  \centering
  \subfloat[Triangulation method]{\label{fig:example-triangle}\includegraphics[width=0.3\textwidth]{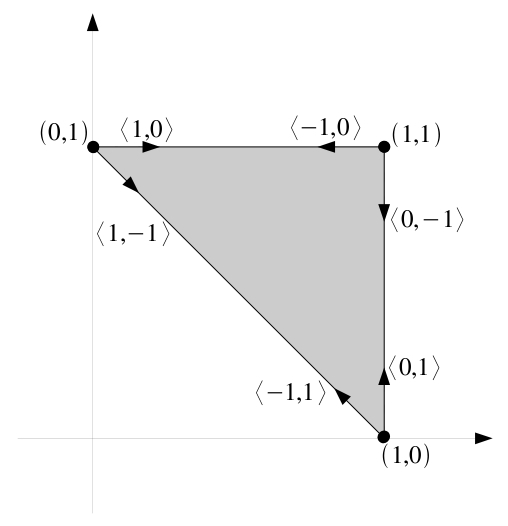}}                
  \subfloat[Cone decomposition method]{\label{fig:example-square}\includegraphics[width=0.3\textwidth]{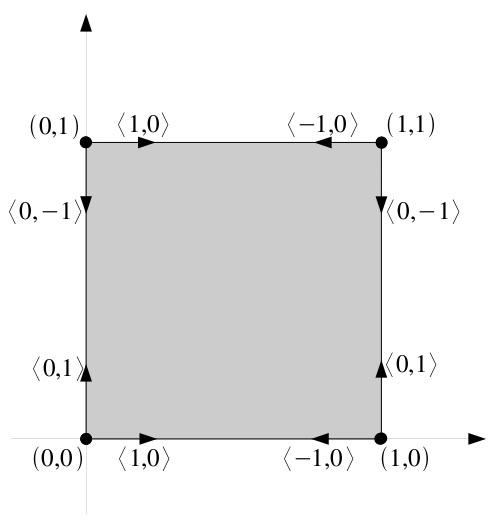}}
  \caption{Example polytopes}
  \label{fig:example-polytopes}
\end{figure}

Take the problem of integrating the polynomial $x+y$ over the triangle $\Delta$ with vertices $s_1=(1,1), s_2=(0,1),$ and $s_3=(1,0)$ in Figure~\ref{fig:example-triangle}. 

The polynomial is already a power of a linear form, and the polytope is a simplex. Because $\ell = (1,1)$ is not regular (it is perpendicular to the edge spanned by $s_2$ and $s_3$), we have to build the index set $K$. Note $\langle \ell, s_1 \rangle = 2, \langle \ell, s_2 \rangle = 1,$ and  $\langle \ell, s_3 \rangle = 1$; pick $K = \{1,2\}$ with $m_1 = 1, m_2 = 2$. We now need to compute two values:

Vertex $s_1$: We are not dividing by zero, we can simply plug vectors into Corollary~\ref{th:brion}, $\frac{\langle \ell, s_1 \rangle ^ 3}{\langle \ell, s_1  - s_2 \rangle ^ 2} = 8.$

Vertex $s_2$: Here, we need to compute a residue. 

\begin{displaymath}
 \Res_{\epsilon=0}\frac{(\epsilon + \langle \ell, s_2 \rangle)^{1+2}}{\epsilon^2 (\epsilon + \langle \ell, s_2 - s_1 \rangle)} = \Res_{\epsilon=0}\frac{(\epsilon +1)^{1+2}}{\epsilon^2 (\epsilon -1)} = -4.
\end{displaymath}

Finally, $\int_{\Delta} (x+y) \d{x}\d{y} = 2! \times \frac{1}{2} \times \frac{1!}{3!} ( 8 - 4) = 2/3$.

\subsubsection{Using the cone decomposition method}
Next, integrate the polynomial $x$ over the unit square in Figure~\ref{fig:example-square} using the cone decomposition algorithm. The polynomial is already a power of a linear form so $\ell = (1,0)$. The polytope has four vertices that we need to consider, and each tangent cone is already simplicial.

Vertex $s_1 = (0,0)$: Because $\langle \ell, s_1 \rangle ^{1+2} = 0,$ the integral on this cone is zero.

Vertex $s_2 = (0,1)$: For the same reason as $s_1$, the integral on this cone is zero.

Vertex $s_3 = (1,0)$: At this vertex, the rays are $u_1 = (0,1), u_2=(-1,0)$. Because $\langle \ell, u_1 \rangle = 0$, we need a perturbation vector $\hat \epsilon$ so that when $\ell := \ell + \hat \epsilon$, we do not divide by zero on any cone (we have to check this cone and the next one). Pick $\hat \epsilon = (\epsilon, \epsilon)$. Then the integral on this cone is
	
\begin{displaymath}
 \frac{M!}{(M+d)!} \vol(\Pi_C) \Res_{\epsilon=0}\frac{(1 + \epsilon)^{1+2}}{\epsilon (-\epsilon) (1 + \epsilon)} = \frac{1!}{(1+2)!} \times 1 \times -2 = -2/6.
\end{displaymath}

Vertex $s_4 = (1,1)$: The rays are $u_1 = (-1,0), u_2= (0,-1)$. Again, we divide by zero, so we perturbate $\ell$ by the same $\hat \epsilon$. The integral on this cone  is
\begin{displaymath}
 \frac{M!}{(M+d)!} \vol(\Pi_C) \Res_{\epsilon=0}\frac{(1 + 2\epsilon)^{1+2}}{\epsilon (\epsilon) (1 + \epsilon)} = \frac{1!}{(1+2)!} \times 1 \times 5 = 5/6.
\end{displaymath}

The integral $\int_P x  \d{x}\d{y} = 0 + 0 -2/6 + 5/6 = 1/2$ as it should be.

\section{How the software works}
\label{section-software}

\latte was originally developed in 2001 as software to study lattice
points of convex polytopes~\cite{latte1}. The algorithms used
combinations of geometric and symbolic computation. Two key
data structures are rational generating functions and cone decompositions,
and it was the first ever implementation of Barvinok's algorithm.
\latte was improved in 2007 with various software and theoretical
modifications, which increased speed dramatically.  This version was
released under the name \texttt{LattE macchiato}; see
\cite{koeppe:irrational-barvinok}.  Now in 2011, our new release \latteInt has extended
its capabilities to include the computation of exact integrals of
polynomial functions over convex polyhedra.  The new integration
functions are C++ implementations of the algorithms provided in
\cite{baldoni-berline-deloera-koeppe-vergne:integration} with
additional technical improvements (including an important new set of
data structures for the manipulation of truncated series). 
A~key distinction between \latteInt and other software tools is that our
algorithms give the exact evaluation of the integral since our
implementation uses exact rational arithmetic.  The code of this
software is freely available at \cite{latteintegrale}

The new implementation of \latteInt allows us to calculate the
integral of a sum of powers of linear forms over an arbitrary
polytope.  Alternatively, we can calculate the integral of a sum of
monomials by decomposing each monomial into a sum of powers of linear
forms using Formula \eqref{eq:decomp-powerlinform}, then integrating
these powers of linear forms.

This section starts with a discussion of our new data structure for
manipulating polynomials and linear forms, then we describe the
format  \latteInt expects for the input polytopes, and we end with a
detailed explanation of the two main algorithms.

\subsection{Input format and data structures}
The input format for the polynomials is identical to that of the \maple
programs of \cite{baldoni-berline-deloera-koeppe-vergne:integration}:
\begin{itemize}
\item A polynomial is represented as a list of its monomials in the form
  \begin{quote}
    \texttt{[{\itshape monomial}$_1$,{\itshape monomial}$_2$,$\dots$]},
  \end{quote}
  where \texttt{{\itshape monomial}$_i$} is represented by
  \begin{quote}
    \texttt{[{\itshape coefficient},[{\itshape exponent-vector}]]}.
  \end{quote}
  For example, $3x_0^2x_1^4x_2^6 + 7x_1^3x_2^5$ is input as \texttt{[[3,[2,4,6]], [7,[0,3,5]]]}. 

\item To deal directly with powers of linear forms, the input format is 
  \begin{quote}
    \texttt{[{\itshape linear-term$_{\hskip 3pt 1}$}, {\itshape linear-term$_{\hskip 3pt 2}$}, $\dots$]},
    \end{quote}
    where \texttt{{\itshape linear-term}$_{\hskip 3pt i}$} is represented by
    \begin{quote}
      \texttt{[{\itshape coefficient},[{\itshape power},[{\itshape coefficient-vector}]]]}.
    \end{quote}
    For example, $3(2x_0 + 4x_1 + 6x_2)^{10} + 7(3x_1 + 5x_2)^{12}$ is input
    as \texttt{[[3,[10,[2,4,6]]], [7,[12,[0,3,5]]]]}. 
\end{itemize}


In \cite{baldoni-berline-deloera-koeppe-vergne:integration}, the integration over 
simplices was first implemented in \maple, and so there was no control over the 
data structures used to store the data. We have implemented the simplex integration 
algorithm in C++ with a sophisticated data structure and have developed a new 
algorithm that integrates over the tangent cones of a polytope. Currently, we are using 
\emph{burst tries}, a data structure designed to have cache-efficient storage and search, 
due to the fact that they are prefix trees with sorted arrays of stored elements as 
leaves~\cite{trip-burst-tries-gastineaau-2006}. Such a data structure is performance-critical 
when computing residues, as a comparison with a linked-list implementation showed. 
In our implementation, each node corresponds to a particular 
dimension or variable and contains the maximal and minimal values of the 
exponent on this dimension.  The node either points to another node a level deeper 
in the tree or a list of sorted elements. 

The input rational polytope $P$ could be given to \latteInt by an h-representation
or v-representation. The input format is the same as in previous versions
and it is explained in \cite{latteintegrale}.  Although the theory we presented earlier works for
both full-dimensional and non-full dimensional rational polytopes, the current release of  \latteInt 
is only guaranteed to do integration and volume computation in full-dimensional polyhedra.
It is worth stressing the old capabilities for counting lattice points still work for polytopes of
all dimensions and  we impose no arbitrary limit on the size or dimension of the input. 
\latteInt relies on {\tt Cddlib} \cite{cddlib-094a} or \fourtitwo \cite{4ti2} for all basic polyhedral 
calculations such as computation of dimension.

Our data structures are specialized for polytopes with vertices of integer coordinates. 
In order to integrate over rational polytopes, we first dilate them and perform a change of variables.
If $P$ is a $d$-dimensional rational polytope and $\alpha P$ is a dilation by $\alpha > 0$ that 
makes $P$ integer, then our software operates on the vertices of $\alpha P$ and rescales
the final integral by the following well-known change of variables:
$$\int_P x_i^{m_i} \ldots x_d^{m_d} \d{m} = \dfrac{1}{\alpha^{d}} \int_{\alpha
  P} \dfrac{1}{\alpha^{m_1+\ldots+m_d}} x_1^{m_1} \ldots x_d^{m_d} \d{m}.$$ 
After this transformation, we apply Formula~\eqref{eq:decomp-powerlinform} to transform
the polynomial into powers of linear forms.

When integrating polytopes other than simplices, there are two options based on the formulas 
presented in Section \ref{preliminar}: 
(i) Triangulate the polytope and apply the algorithm for each simplices individually, or 
(ii) Triangulate each tangent cone and integrate each one using the cone decomposition algorithm. 
Therefore,  a key step in all our computations is to find either a triangulation of the polytope or a triangulation 
of each of its tangent cones. Once more, this step relies on {\tt Cddlib} or \fourtitwo, because when we triangulate
we compute a regular triangulation via a convex hull \cite{latte1,koeppe:irrational-barvinok}. 
We now explore the two integration algorithms in greater detail.

%
	          

\subsection{Integrating powers of linear forms by polytope triangulation}

After we decompose the polynomial to a sum of powers of linear forms and after finding a 
triangulation of the polytope,  Algorithm~\ref{alg:integratePolytopeTriangulation} loops over 
these two sets and integrates each linear form/simplex pair individually using Corollaries \ref{th:brion} and~\ref{th:residue}.

\begin{algorithm}                      
\caption{ Integrate using polytope triangulation}
\label{alg:integratePolytopeTriangulation}
\begin{algorithmic}                    
\REQUIRE $F = \sum c_j \langle \ell_j, x \rangle^{M_j}$, sum of powers of linear forms.
\REQUIRE $P$, a full-dimensional polytope.
\ENSURE integral of the linear forms $F$ over the polytope P.
\STATE $\mathit{integral} \gets 0$ \COMMENT{ $\mathit{integral}$ is a rational data type}
\STATE let $T$ be a list of simplices that form a triangulation of $P$
\FORALL{simplices $\Delta$ in $T$}
	\FORALL{linear forms $c \langle \ell, x \rangle^M$ in $F$}
		\IF{$\ell$ is regular on $\Delta$}
			\STATE $\mathit{integral} \gets \mathit{integral} + c \times \mathit{integrateSimplexRegular}(\ell, M, \Delta)$
		\ELSE
			\STATE $\mathit{integral} \gets \mathit{integral} + c \times \mathit{integrateSimplexResidue}(\ell, M, \Delta)$
		\ENDIF
	\ENDFOR
\ENDFOR
\RETURN $\mathit{integral}$
\end{algorithmic}
\end{algorithm}

In Algorithm~\ref{alg:integratePolytopeTriangulation}, the  linear forms are 
represented as a burst trie, the triangulations are stored in a linked list, and 
each simplex is a simple two-dimensional array containing the vertices $s_1, \dots, s_{d+1}$.

When $\ell$ is regular on $\Delta$, the \emph{integrateSimplexRegular} function plugs 
in numbers and vectors into Corollary~\ref{th:brion}. Also, the terms in the numerator and 
denominator are evaluated in a rational data type, and so no floating-point divisions are performed.

When $\ell$ is not regular, the \emph{integrateSimplexResidue} function (Algorithm~\ref{alg:integrateSimplexResidue}) 
applies Corollary~\ref{th:residue} and must find an index set $K \subset \{1,\dots,d+1\}$ of different poles 
$t=1/\langle \ell,s_k\rangle $, and compute $|K|$ residues. Let $k \in K$  and let $m_k$ denote the order of the pole, i.e.,

\begin{center}
$m_k=\#\{\, i\in \{1,\dots ,d+1\}:\langle \ell,s_i\rangle =\langle \ell,s_k\rangle \, \}$.
\end{center}

The problem has now been reduced to evaluating Formula~\eqref{eq:simplexResidueCoeff}.

\begin{equation}
\Res_{\epsilon=0}\frac{(\epsilon+\langle \ell,s_k\rangle )^{M+d}}{\epsilon ^{m_k}\prod\limits_{\substack{i\in K \\i\neq k}} (\epsilon+\langle \ell,s_k-s_i\rangle )^{m_i}}=[\epsilon^{m_k-1}]\frac{(\epsilon+\langle \ell,s_k\rangle )^{M+d}}{\prod\limits_{\substack{i\in K \\i\neq k}} (\epsilon+\langle \ell,s_k-s_i\rangle)^{m_i}},
\label{eq:simplexResidueCoeff}
\end{equation}
where $[\epsilon^{a}]p$ means the coefficient of $\epsilon^a$ in the Laurent series of expression $p$.

To compute Formula~\eqref{eq:simplexResidueCoeff}, we expand the polynomial in terms of $\epsilon$ in the numerator truncated to degree $m_k -1$. We then find the first $m_k$ terms in the polynomial expansion of each $1/ (\epsilon+\langle \ell,s_k - s_i\rangle)^{m_i}, i \neq k$ term using the general binomial theorem. To make the notation easy, let $b = \langle\ell, s_k - s_i\rangle \in \Z$, then the degree $m_k-1$ polynomial of $(\epsilon+b)^{-m_k}$ is
\begin{displaymath}
\label{factorPolynomial}
\hat{p}(\epsilon) = \alpha_0 \epsilon^0 b^{-m_i} + \alpha_1 \epsilon^1 b^{-m_i-1} + \dots + \alpha_{m_k-1} \epsilon^{m_k-1} b^{-m_i-m_k+1},\; \alpha_j = \binom{m_i + j - 1}{m_i - 1}(-1)^j.
\end{displaymath}

This is a polynomial in $\epsilon$ with rational coefficients. For efficiency reasons, we factor $\hat{p}(\epsilon) = \frac{1}{b^{m_i + m_k -1}}p(\epsilon),\,  p \in \Z[\epsilon]$. The \emph{integrateSimplexResidue} and \emph{truncatedMultiply} functions both implement these ideas.

\begin{algorithm}                      
\caption{integrateSimplexResidue}
\label{alg:integrateSimplexResidue}
\begin{algorithmic}                    
\REQUIRE $\ell$, coefficients of a linear form and $M$, integer power.
\REQUIRE $\Delta$, simplex with vertices $s_1, \dots, s_{d+1}$.
\ENSURE The integral of $\langle \ell, x\rangle^{M}$ over $\Delta$.

\STATE Let $p_1 \gets 1, p_2 \gets 1$ be polynomials in $\epsilon$.
\STATE Let $\mathit{rf} \gets 1$ be a rational data type.
\STATE Let $\mathit{sumResidue} \gets 0$ be a rational data type.
\STATE Make the index set $K$ of unique poles.

\FORALL{$k$ in $K$}
	\STATE $\mathit{rf} \gets 1$
	\STATE $p_1 \gets$ the expansion of $(\epsilon + \langle\ell, s_k\rangle)^{M+d}$ up to degree $m_k-1$
	\COMMENT{$p_1 \in \Z[\epsilon]$}
	
	\FORALL{$i$ in $K$ and $i \neq k$}
		\STATE $\mathit{rf} \gets \mathit{rf} \times  \langle\ell, s_k\rangle^{- (m_i +m_k -1)}$
		\STATE $p_2 \gets$ the expansion of $(\epsilon + \langle\ell, s_k - s_i\rangle)^{-m_i}$ up to degree $m_k -1$ with $\langle\ell, s_k\rangle^{- (m_i +m_k -1)}$ factored out.
		\COMMENT{$p_2 \in \Z[\epsilon]$}

		\STATE $p_1 \gets \mathit{truncatedMultiply}(p_1, p_2, m_k -1$)

	\ENDFOR 
	\STATE Let $c$ be the coefficient of the degree $m_k-1$ term in $p_1(\epsilon)$.
	\STATE $\mathit{sumResidue} \gets \mathit{sumResidue} + \mathit{rf} \times c$ 

\ENDFOR 
\RETURN $abs(det(s_1 - s_{d+1}, \dots, s_d - s_{d+1})) \times \frac{M!}{(M+d)!} \times \mathit{sumResidue}$
\end{algorithmic}
\end{algorithm}

Finally, Algorithm \emph{truncatedMultiply(p,q,k)} takes two polynomials $p, q \in \Z[\epsilon]$ and 
returns their  product up to and including terms of degree $k$. Our implementation is very simple (e.g., not using any
special multiplication algorithms) but the cache-efficient use of the burst tries leads to speed-ups wheatn compared to a
naive implementation with arrays. We note that asymptotically faster multiplication algorithms exists (see, e.g., \cite{Guillaumeetal}),
which might lead to further improvements.

\subsection{Integrating powers of linear forms by cone decomposition}

After triangulating each tangent cone into simplicial cones, the computation is very similar to the 
polytope-triangulation case: if $\ell$ is regular on the rays of the cone, we plug in values into 
Corollary~\ref{oneconeoneplf}, else we perturb $\ell$ and perform a residue calculation. 
Algorithm~\ref{alg:integratePolytopeTangentCone} implements this idea.

\begin{algorithm}                      
\caption{Integrate using the cone decomposition method}
\label{alg:integratePolytopeTangentCone}
\begin{algorithmic}                    
\REQUIRE $F = \sum c_j \langle \ell_j, x \rangle^{M_j}$, powers of linear forms.
\REQUIRE $P$, a full-dimensional polytope.
\ENSURE integral of the linear forms over $P$.
\STATE $\mathit{integral} \gets 0$ \COMMENT{ $\mathit{integral}$ is a rational data type}
\STATE Let $C$ be a list of tangent cones $P$.
\STATE Make $T$ be a list of triangulated cones in $C$.
\COMMENT{A cone in $T$ is in the form $(s; u_1, \dots, u_d)$, where $s$ is a vertex and $u_i$ are rays}

\FORALL{linear forms $c \langle \ell, x \rangle^M$ in $F$}
	\STATE Let $R \subseteq T$ be cones that $\ell$ is regular on.

	\FORALL{$ (s; u_1, \dots, u_d)$ in $R$}
		\STATE $\mathit{integral} \gets \mathit{integral} + c \times \mathit{integrateConeRegular}(\ell, M, s, u_1, \dots, u_d)$
	\ENDFOR

	\STATE Pick $\hat \epsilon = \epsilon(a_1, \dots, a_d)$ where $a_i \in \Z$ so that we do not divide by zero.

	\FORALL{$ (s; u_1, \dots, u_d)$ in $T\setminus R$}
		\STATE $\mathit{integral} \gets \mathit{integral} + c \times \mathit{integrateConeResidue}(\ell, M, \hat \epsilon, s, u_1, \dots, u_d)$
	\ENDFOR

\ENDFOR
\RETURN $\mathit{integral}$
\end{algorithmic}
\end{algorithm}

In \emph{integrateConeResidue}, $\ell$ is perturbed by setting $\ell:=\ell + \hat \epsilon$, where $\hat \epsilon$ 
is a vector in terms of $\epsilon$ with coefficients picked on the moment curve with alternating signs. We repeatedly 
pick a random  $t \in \Z_+$ and set $\hat \epsilon_i = t^{i-1}(-1)^{i-1} \epsilon$ for $i = 1,2,\dots,d$ until 
$\langle -(\ell + \hat \epsilon), u_i \rangle$ is non-zero for every simple cone at every vertex in Corollary~\ref{oneconeoneplf}. 
Then the residues are computed using the general binomial theorem and truncated series multiplication like in \emph{integrateSimplexResidue}.

\subsection{A special case: computing volumes}
Computing the volume of a polytope is equivalent to integrating the monomial $1$ over 
the polytope. We again have the two same options when computing volumes as we did 
when computing integrals. Instead of using the algorithms above, we can simplify the computation.
In the triangulation based approach, we find a triangulation of the polytope and sum 
the volume of each simplex. The volume of a specific simplex is obtained by taking a determinant. 
In the cone decomposition approach, we triangulate each tangent cone and apply 
Corollary~\ref{oneconeoneplf} with $M=0$ and any random vector $\ell$. If we do 
divide by zero, instead of finding residues, we simply pick a new random $\ell$ and 
start the computation over.

\section{Experiments}
\label{section-integration-experiments}

We did thorough testing of the implementation and performed new computational benchmarks. We report on four different test classes:

\begin{enumerate}
	\item We expand the computational limits for integrating over simplices described in \cite{baldoni-berline-deloera-koeppe-vergne:integration}. 
	\item Next, we integrate random monomials over three families of polytopes: (1) simple, (2) simplicial, and (3) neither simple nor simplicial.
	\item Because our volume methods are optimized versions of the integration methods, we also compute the volumes of the same polytopes in the last case above.
	\item Finally, we compare \latteInt to other software tools and computational studies \cite{bueler-enge-fukuda-2000:exact-volume,CUBPACK, korenblit-2006:integration-over-convex-polyhedron,Ong-2003:finding-exact-volume-polyhedron}.
\end{enumerate}

Our integration and volume experiments input data, along with running times and the results of integration and volumes are available on the 
\latte website \cite{latteintegrale}.

\subsection{Integration over simplices}
In \cite{baldoni-berline-deloera-koeppe-vergne:integration}, the theory of integration over simplices was developed and
a fair amount of \maple experiments were carried out to show the potential of the methods. In this section, the experiments we performed
clearly indicate that this C++ implementation is at least two orders of magnitude faster than the preceding \maple code; 
compare Tables 5 and 6 in \cite{baldoni-berline-deloera-koeppe-vergne:integration} with 
Tables~\ref{table:simplex-linear-form} and~\ref{table:simple-polynomial} in this paper. In Table~\ref{table:simplex-linear-form}, 
we used \maple to generate powers of random linear forms and randomly generated simplices. The coefficients of each linear form were 
picked uniformly over $[0,100] \cap \Z$. We did the integration using \latteInt. 

\begin{table}
\centering
\caption{Average integration time plus or minus one standard deviation when integrating one power of a linear form over a random $d$-simplex (in seconds over 50 random forms)}
\label{table:simplex-linear-form}
\small
\begin{tabular}{llllllll}
\toprule
   & \multicolumn{7}{c}{Exponent $M$} \\
\cmidrule(r){2-8}	
$d$  & \multicolumn{1}{c}{2} & \multicolumn{1}{c}{10} & \multicolumn{1}{c}{20} & \multicolumn{1}{c}{50} & \multicolumn{1}{c}{100} & \multicolumn{1}{c}{300} & \multicolumn{1}{c}{1000} \\ \hline

 10  & $0.00\pm0.00$ & $0.00\pm0.00$ & $0.00\pm0.00$ & $0.00\pm0.00$ & $0.00\pm0.00$ & $0.00\pm0.00$ & $0.01\pm0.00$\\
 20  & $0.00\pm0.00$ & $0.00\pm0.00$ & $0.00\pm0.00$ & $0.00\pm0.00$ & $0.00\pm0.00$ & $0.00\pm0.00$ & $0.01\pm0.00$\\
 50  & $0.00\pm0.00$ & $0.00\pm0.00$ & $0.00\pm0.00$ & $0.00\pm0.00$ & $0.00\pm0.01$ & $0.01\pm0.00$ & $0.03\pm0.01$\\
 100  & $0.02\pm0.00$ & $0.02\pm0.00$ & $0.02\pm0.00$ & $0.02\pm0.00$ & $0.02\pm0.00$ & $0.03\pm0.00$ & $0.09\pm0.00$\\
 300  & $0.35\pm0.01$ & $0.36\pm0.01$ & $0.36\pm0.01$ & $0.36\pm0.01$ & $0.38\pm0.01$ & $0.42\pm0.01$ & $0.66\pm0.02$\\
 400  & $0.78\pm0.02$ & $0.79\pm0.03$ & $0.79\pm0.02$ & $0.80\pm0.03$ & $0.82\pm0.03$ & $0.90\pm0.03$ & $1.25\pm0.04$\\
\bottomrule
\end{tabular}
\end{table}

Next, in Table~\ref{table:simple-polynomial}, we used \maple to generate monomials and simplices. We again integrate using
\latteInt. We measure time from the start of program execution to termination, which includes file I/O, system calls, child process time, 
the time to find tangent-cones, and triangulation time. All triangulations were computed with the software package \texttt{cddlib} version 0.94f \cite{cddlib-094a}. 
All computations  were performed on a 64-bit Ubuntu machine with 64 GB of RAM and eight Dual Core AMD Opteron 880 processors.
We applied a 600-second maximum running time to this program; tasks taking longer are not benchmarked.

\begin{table}
\centering
\caption{Average and standard deviation of integration time in seconds of a random monomial of prescribed degree by decomposition into linear forms over a $d$-simplex (average over 50 random forms)}
\label{table:simple-polynomial}

\small
\begin{tabular}{lrrrrrrrrrrr}
\toprule
& \multicolumn{11}{c}{Degree} \\ 
\cmidrule(r){2-12}	
$d$ & \multicolumn{1}{c}{1} & \multicolumn{1}{c}{2} & \multicolumn{1}{c}{5} & \multicolumn{1}{c}{10} & \multicolumn{1}{c}{20} & \multicolumn{1}{c}{30} & \multicolumn{1}{c}{40} & \multicolumn{1}{c}{50} & \multicolumn{1}{c}{100} & \multicolumn{1}{c}{200} & \multicolumn{1}{c}{300} \\ \hline
 2  & 0.0 & 0.0 & 0.0 & 0.0 & 0.0 & 0.0 & 0.0 & 0.0 & 0.1 & 1.0 & 3.8\\ 
  & 0.0 & 0.0 & 0.0 & 0.0 & 0.0 & 0.0 & 0.0 & 0.0 & 0.1 & 0.4 & 1.7\\
 3  & 0.0 & 0.0 & 0.0 & 0.0 & 0.0 & 0.0 & 0.1 & 0.2 & 2.3 & 38.7 & 162.0\\ 
  & 0.0 & 0.0 & 0.0 & 0.0 & 0.0 & 0.0 & 0.0 & 0.1 & 1.4 & 24.2 & 130.7\\
 4  & 0.0 & 0.0 & 0.0 & 0.0 & 0.0 & 0.1 & 0.4 & 0.7 & 22.1 & -- & --\\ 
  & 0.0 & 0.0 & 0.0 & 0.0 & 0.0 & 0.1 & 0.3 & 0.7 & 16.7 & -- & --\\
 5  & 0.0 & 0.0 & 0.0 & 0.0 & 0.1 & 0.3 & 1.6 & 4.4 & -- & -- & --\\ 
  & 0.0 & 0.0 & 0.0 & 0.0 & 0.0 & 0.2 & 1.3 & 3.5 & -- & -- & --\\
 6  & 0.0 & 0.0 & 0.0 & 0.0 & 0.1 & 1.1 & 4.7 & 15.6 & -- & -- & --\\ 
  & 0.0 & 0.0 & 0.0 & 0.0 & 0.1 & 1.0 & 4.3 & 14.2 & -- & -- & --\\
 7  & 0.0 & 0.0 & 0.0 & 0.0 & 0.2 & 2.2 & 12.3 & 63.2 & -- & -- & --\\ 
  & 0.0 & 0.0 & 0.0 & 0.0 & 0.2 & 1.7 & 12.6 & 66.9 & -- & -- & --\\
 8  & 0.0 & 0.0 & 0.0 & 0.0 & 0.4 & 4.2 & 30.6 & 141.4 & -- & -- & --\\ 
  & 0.0 & 0.0 & 0.0 & 0.0 & 0.3 & 3.0 & 31.8 & 127.6 & -- & -- & --\\
 10  & 0.0 & 0.0 & 0.0 & 0.0 & 1.3 & 19.6 & -- & -- & -- & -- & --\\ 
  & 0.0 & 0.0 & 0.0 & 0.0 & 1.4 & 19.4 & -- & -- & -- & -- & --\\
 15  & 0.0 & 0.0 & 0.0 & 0.1 & 5.7 & -- & -- & -- & -- & -- & --\\ 
  & 0.0 & 0.0 & 0.0 & 0.0 & 3.6 & -- & -- & -- & -- & -- & --\\
 20  & 0.0 & 0.0 & 0.0 & 0.2 & 23.3 & -- & -- & -- & -- & -- & --\\ 
  & 0.0 & 0.0 & 0.0 & 1.3 & 164.8 & -- & -- & -- & -- & -- & --\\  
 30  & 0.0 & 0.0 & 0.0 & 0.6 & 110.2 & -- & -- & -- & -- & -- & --\\ 
  & 0.0 & 0.0 & 0.1 & 4.0 & 779.1 & -- & -- & -- & -- & -- & --\\
 40  & 0.0 & 0.0 & 0.0 & 1.0 & -- & -- & -- & -- & -- & -- & --\\ 
  & 0.0 & 0.0 & 0.3 & 7.0 & -- & -- & -- & -- & -- & -- & --\\
 50  & 0.0 & 0.0 & 0.1 & 1.8 & -- & -- & -- & -- & -- & -- & --\\ 
  & 0.0 & 0.1 & 0.5 & 12.9 & -- & -- & -- & -- & -- & -- & --\\

\bottomrule
\end{tabular}
\end{table}




\subsection{Integration over general polytopes}

We tested the triangulation and cone decomposition integration methods on polytopes and their duals across dimension, vertex counts, and over monomials with different degrees.
 For each polytope dimension and vertex count we constructed 50 random polytopes by taking the convex hull of random points using {\tt Polymake} \cite{polymake-software}. 
 For primal polytopes of dimension $d$, the number of vertices considered goes from $d+2$ to $d+25$.
 When zero is not in the interior of the polytope, we translated the centroid to the origin before constructing the dual polytope. Then we integrated each polytope and its corresponding 
 dual polytope over a new random monomial of a set degree. Because of the construction method, most primal polytopes are simplicial and the duals are mostly simple polytopes. 
 We also integrated over G. M. Ziegler's database of polytopes  \cite{ziegler-db-website}, which contains polytopes that are not simplicial nor simple.

\begin{figure}
  \centering  
  \subfloat{\label{fig:i3p}\includegraphics[width=1.0\textwidth]{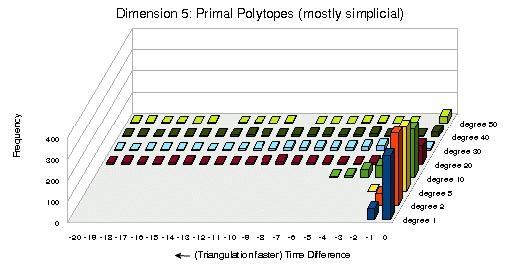}}  
  \\              
  \subfloat{\label{fig:i3d}\includegraphics[width=1.0\textwidth]{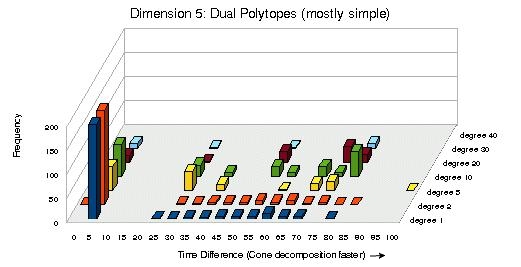}}
  \caption{Histogram of the time difference between the triangulation and cone-decomposition methods for integrating  over random polytopes in dimension 5}
  \label{fig:integration-random-graphs1}
\end{figure}

\begin{figure}
  \centering  
  \subfloat{\label{fig:ip}\includegraphics[width=1.0\textwidth]{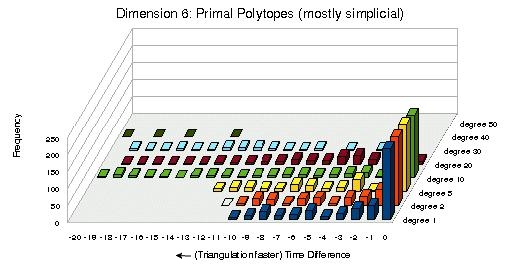}}
  \\                
  \subfloat{\label{fig:id}\includegraphics[width=1.0\textwidth]{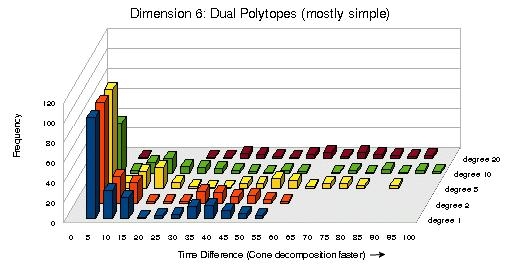}}
  \caption{Histogram of the time difference between the triangulation and cone-decomposition methods for integrating  over random polytopes in dimension 6}
  \label{fig:integration-random-graphs3}
\end{figure}

\begin{figure}
  \centering  
  \subfloat{\label{fig:ip7}\includegraphics[width=1.0\textwidth]{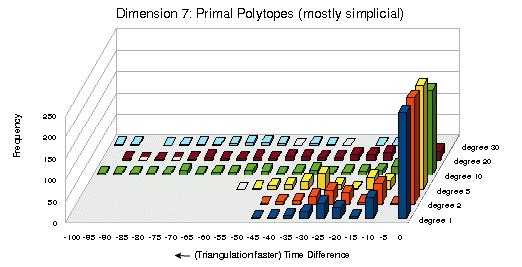}} 
  \\               
  \subfloat{\label{fig:id7}\includegraphics[width=1.0\textwidth]{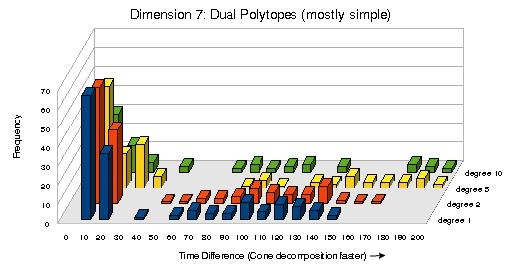}}
  \caption{Histogram of the time difference between the triangulation and cone-decomposition methods for integrating  over random polytopes in dimension 7}
  \label{fig:integration-random-graphs5}
\end{figure}

\begin{figure}
  \centering  
  \subfloat{\label{fig:izp}\includegraphics[width=1.0\textwidth]{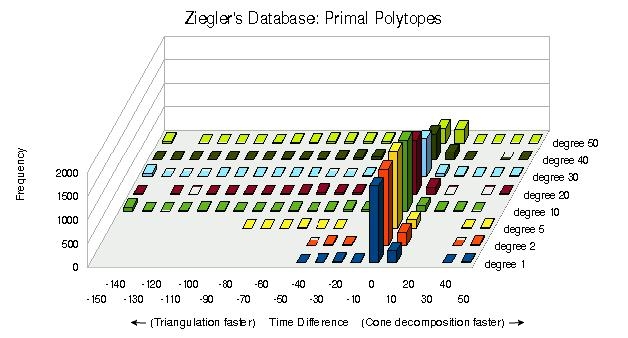}}   
  \\
  \subfloat{\label{fig:izd}\includegraphics[width=1.0\textwidth]{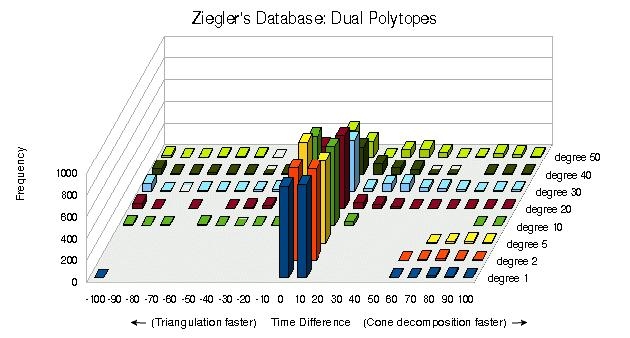}}
  \caption{Histogram of the time difference between the triangulation and cone-decomposition methods for integrating over the polytopes in Ziegler's Database}
  \label{fig:integration-zielger}
\end{figure}

\begin{description}

\item[Simple and simplicial polytopes]
We present the results for dimensions 5, 6, and 7.  We tested both algorithms on the primal polytopes starting with their v-representation. For their duals we tested the triangulation and cone decomposition methods starting from their h-representations.   We also did experiments in dimension 3 and 4 but
the numbers are too close to each other to show a clear trend of which is the fastest method (see \cite{extratables}).
We only report those test classes for which every polytope in the test class finished under 600 seconds for both the triangulation and the cone-decomposition method.

In Figures \ref{fig:integration-random-graphs1}, \ref{fig:integration-random-graphs3}, and \ref{fig:integration-random-graphs5}, 
we display histograms that on the $x$-axis plot the time difference between the two integration methods, the $y$-axis shows the degree
of monomials and the $z$-axis presents the number of random polytopes (in the respective dimensions 5, 6, 7). A particular solid bar in position $(a^*,b^*)$ tallies 
the number of random polytopes for which the time difference between the two algorithms was $a^*$ seconds when integrating a
a monomial of degree $b^*$. We define the \emph{time difference} as the time taken by the triangulation method \emph{minus} the time
taken by the cone-decomposition method. There is one histogram for the primal polytopes and one for the dual polytopes.
For example, in Figure \ref{fig:integration-random-graphs1} on the top (primal), there are eight colors on the bars, one for each degree (corresponds to a row).
The bars in the degree 50 row, show the average time differences are always less than or equal to zero. This shows that the cone method was
slower in the majority of simplicial problems. For comparison, on the bottom (dual), most of the mass has positive time difference, which
indicates the cone method wins over the integration of higher degree monomials. More tables will be available online \cite{extratables}.

In conclusion, our results for dimension higher than four show integrating monomials coincide qualitatively with the observation of 
\cite{bueler-enge-fukuda-2000:exact-volume} for volume computation (polynomial of degree zero): the triangulation 
method is faster for simplicial polytopes (mass on histograms is concentrated on negative time differences)  while the cone 
decomposition is faster for simple polytopes (mass on histograms is concentrated on positive time differences).  The trends are very clear
while in dimension less than four, the timings are too close to each other to give a clear-cut trend.



\item[Zero-one polytopes]
In Figure \ref{fig:integration-zielger}, we present another histogram comparing the cone decomposition and triangulation 
methods on Ziegler's database of polytopes \cite{ziegler-db-website}, which contains many zero-one polytopes and a few 
other polytopes. We translate each polytope so that its centroid is the origin, thus its dual is well-defined. Then for each 
polytope and its dual, we integrate 50 random polynomials of a set degree. We skipped non-full-dimensional polytopes 
and a few others that did not finish within 30 minutes. The figure displays the histogram of the differences in running times 
between the two integration methods for monomials of eight different degrees (1, 2, 5, 10, 20, 30, 40, and 50). As before, 
the $z$-axis is frequency, $x$-axis is the range of the time differences, and $y$-axis gives the degree used. For example, there 
are over 200 cases where integrating a degree  40 monomial over the dual polytopes took the triangulation method 90 
seconds more than the cone-decomposition method.

The behavior  we observed before for simple vs. simplicial polytopes still mostly holds for these tests, except
we see that the two methods finish within 20 seconds of each other most of the time (mass of histogram is centered at 0). 
The variation is then not as strong as before. 

\end{description}

\subsection{Volume Experiments}
\label{section-volume-experiments}

Volume computation is an important special case of integration that has received attention by several researchers, thus 
we also tested the triangulation and cone decomposition methods on the same database of random polytopes and their 
duals, and on Ziegler's database to see the performance of volume evaluation.

\begin{figure}
  \centering  
  \subfloat{\label{fig:rvp}\includegraphics[width=1.0\textwidth]{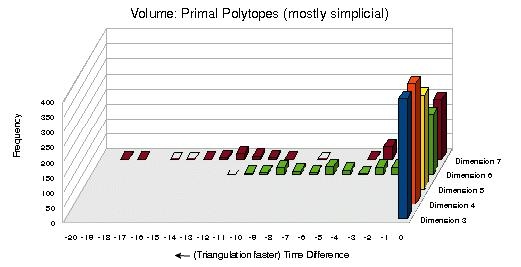}}   
  \\
  \subfloat{\label{fig:rvd}\includegraphics[width=1.0\textwidth]{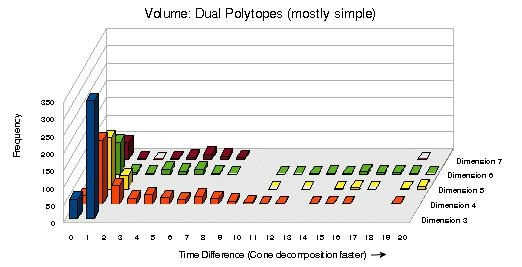}}
  \caption{Histogram of the time difference between the triangulation and cone-decomposition methods for finding the volume of random polytopes}
  \label{fig:random-volume}
\end{figure}

\begin{description}
\item[Simple and simplicial polytopes]

Each test class contains 50 polytopes  for each dimension and we only considered tests where both methods  finished within 600 seconds for the same polytope.
While the triangulation method is still faster for simplicial polytopes and the cone-decomposition method is faster for simple polytopes,
the histograms in Figure \ref{fig:random-volume} suggest that both methods finish quite close to each other in small dimension. When
dimension starts growing there is a more pronounced difference between the methods (i.e., the mass of
the histogram is more spread toward positive or negative values of the time difference).

\item[Zero-one polytopes] In
  Table~\ref{tabel:lawrence-ziegler-volume-1} and
  \ref{tabel:lawrence-ziegler-volume-2}, we apply the triangulation
  and cone decomposition volume methods to Ziegler's database
  \cite{ziegler-db-website} and their duals. If a polytope did not
  contain the origin, we centered it so that its dual is
  defined. Again, we skipped non-full-dimensional polytopes and a few
  others that did not finish within 30 minutes.
Faster timings are shown in bold.
When computing volumes of primal polytopes in Ziegler's database, triangulation is faster more often. But for  finding the volume of the dual polytopes there is no clear faster method.
\end{description}

\begin{table}
\centering
\caption[Volume of Ziegler's Database: Part I]{The triangulation vs.~cone-decomposition method for finding volumes in Ziegler's Database: Part I}

\label{tabel:lawrence-ziegler-volume-1}
\tabcolsep 2.0pt
\small
\begin{tabular}{lrrrrrrr}
\toprule
& & \multicolumn{3}{c}{Primal} &  \multicolumn{3}{c}{Dual} \\
		\cmidrule(r){3-5} \cmidrule(r){6-8}
& & & \multicolumn{2}{c}{Time (sec)} & & \multicolumn{2}{c}{Time (sec)} \\
		\cmidrule(r){4-5} \cmidrule(r){7-8}
Polytope & Dim. & Vertices & Cone. &  Triang. & Vertices & Cone. & Triang. \\
\hline                    
3simp3simp.vrep.latte		&6		&44	&\timeFaster{ 5.61	}&\timeSlower{ 6.10		}&32		&\timeFaster{ 1.11		}&\timeSlower{ 1.15 }\\
cyclic\_4\_8.vrep.latte		&4		&8		&\timeSlower{ 0.09		}&\timeFaster{ 0.06	}&20		&\timeFaster{ 0.02		}&\timeSlower{ 0.10}\\
neighborly\_4\_8.vrep.latte&4		&8		&\timeSlower{ 0.12	}&\timeFaster{ 0.03		}&20		&\timeFaster{ 0.03		}&\timeSlower{ 0.06}\\
SharirCube.vrep.latte		&3		&8		&0.03		& 0.03	&6		&\timeSlower{0.11		}&\timeFaster{ 0.02}\\
HC\_6-32.vrep.latte		&6		&32		&\timeSlower{ 2.29		}&\timeFaster{ 2.06	}&44		&\timeSlower{3.25		}&\timeFaster{ 3.22}\\
HC\_7-64.vrep.latte		&7		&64		&\timeFaster{ 13.42		}&\timeSlower{75.85	}&78		&\timeFaster{ 61.68		}&\timeSlower{ 762.12}\\
HC\_8-128.vrep.latte		&8		&128		&\timeFaster{ 85.85		}&--						&144		&\timeFaster{ 15007.50		}&--\\
MJ\_16-17.vrep.latte		&16		&17	&\timeSlower{ 2.60		}&\timeFaster{ 2.48	}&17		&\timeSlower{ 0.07		}&\timeFaster{ 0.04}\\
OA\_5-10.vrep.latte		&5		&10		&\timeSlower{ 0.22		}&\timeFaster{ 0.08	}&22		&\timeSlower{ 0.18		}&\timeFaster{ 0.11}\\
OA\_5-18.vrep.latte		&5		&18		&\timeSlower{ 0.46		}&\timeFaster{ 0.32	}&19		&\timeSlower{ 0.34		}&\timeFaster{ 0.10}\\
OA\_5-24.vrep.latte		&5		&24		&\timeSlower{ 0.81		}&\timeFaster{ 0.58	}&18		&\timeSlower{ 0.22		}&\timeFaster{ 0.13}\\
OA\_6-13.vrep.latte		&6		&13		&\timeSlower{ 0.53		}&\timeFaster{ 0.20	}&56		&\timeFaster{ 0.52		}&\timeSlower{ 5.37}\\
OA\_7-18.vrep.latte		&7		&18		&\timeSlower{ 3.36		}&\timeFaster{ 0.82	}&146		&\timeFaster{ 13.96		}&\timeSlower{ 1827.83}\\
OA\_8-25.vrep.latte		&8		&25		&\timeSlower{ 38.55	}&\timeFaster{ 10.44		}&524		&\timeFaster{ 4116.93	}&--\\
OA\_9-33.vrep.latte		&9		&33		&--		&\timeFaster{ 648.77		}&1870		&--		&--\\
AS\_6-18.vrep.latte		&6		&18		&\timeSlower{ 1.26		}&\timeFaster{ 0.52	}&121		&\timeFaster{ 1.40		}&\timeSlower{ 65.63}\\
BIR3\_4-6.vrep.latte		&4		&6			&\timeSlower{ 0.12		}&\timeFaster{ 0.02	}&9		& 0.01		&0.01\\
BIR4\_9-24.vrep.latte		&9		&24	&\timeSlower{ 6.26		}&\timeFaster{ 2.22	}&16		&\timeSlower{ 1.42		}&\timeFaster{ 0.17}\\
BIR5\_16-120.vrep.latte		&16	&120	&--							 &--						&25		&--							&\timeFaster{ 488.78}\\

\bottomrule

\end{tabular}
\end{table}

\begin{table}
\centering
\caption[Volume of Ziegler's Database: Part I]{The triangulation vs.~cone-decomposition method for finding volumes in Ziegler's Database: Part II}

\label{tabel:lawrence-ziegler-volume-2}
\tabcolsep 2.0pt
\small
\begin{tabular}{lrrrrrrr}
\toprule
& & \multicolumn{3}{c}{Primal} &  \multicolumn{3}{c}{Dual} \\
		\cmidrule(r){3-5} \cmidrule(r){6-8}
& & & \multicolumn{2}{c}{Time (sec)} & & \multicolumn{2}{c}{Time (sec)} \\
		\cmidrule(r){4-5} \cmidrule(r){7-8}
Polytope & Dim. & Vertices & Cone. &  Triang. & Vertices & Cone. & Triang. \\
\hline                    

CF\_10-11.vrep.latte		&10	&11		&\timeSlower{ 0.37		}&\timeFaster{ 0.33	}&11		&\timeSlower{ 0.03		}&\timeFaster{ 0.00}\\
CF\_4-5.vrep.latte		&4		&5		&\timeSlower{ 0.01		}&\timeFaster{ 0.02		}&5		&0.02		&0.02\\
CF\_5-6.vrep.latte		&5		&6		&\timeSlower{ 0.04		}&\timeFaster{ 0.03		}&6		&0.00		&0.00\\
CF\_6-7.vrep.latte		&6		&7		&\timeFaster{ 0.04		}&\timeSlower{ 0.05		}&7		&\timeSlower{ 0.01		}&\timeFaster{ 0.00}\\
CF\_7-8.vrep.latte		&7		&8		&\timeFaster{ 0.08		}&\timeSlower{ 0.09		}&8		&\timeSlower{ 0.02		}&0.02\\
CF\_8-9.vrep.latte		&8		&9		&\timeSlower{ 0.14		}&\timeFaster{ 0.13		}&9		&\timeSlower{ 0.02		}&\timeFaster{ 0.01}\\
CF\_9-10.vrep.latte		&9		&10	&\timeSlower{ 0.22		}&\timeFaster{ 0.20		}&10		&\timeFaster{ 0.01		}&0.02\\
CRO\_3-6.vrep.latte		&3		&6		&\timeSlower{ 0.04		}&\timeFaster{ 0.01		}&8		&\timeFaster{ 0.00		}&\timeSlower{ 0.02}\\
CRO\_4-8.vrep.latte		&4		&8		&\timeSlower{ 0.12		}&\timeFaster{ 0.05		}&16		&\timeFaster{ 0.00		}&\timeSlower{ 0.03}\\
CRO\_5-10.vrep.latte		&5		&10	&\timeSlower{ 0.17		}&\timeFaster{ 0.10		}&32		&\timeFaster{ 0.01		}&\timeSlower{ 0.33}\\
CUT3\_3-4.vrep.latte		&3		&4		&0.00							 &\timeSlower{ 0.01		}&4		&0.01		&0.01\\
CUT4\_6-8.vrep.latte		&6		&8		&\timeSlower{ 0.16		}&\timeFaster{ 0.06		}&16		&\timeFaster{ 0.00		}&\timeSlower{ 0.05}\\
CUT5\_10-16.vrep.latte	&10	&16	&\timeSlower{ 2.72		}&\timeFaster{ 0.94		}&56		&\timeFaster{ 38.38		}&\timeSlower{ 2046.74}\\
CYC\_5-8.vrep.latte		&5		&8		&\timeSlower{ 0.11		}&\timeFaster{ 0.05		}&20		&\timeFaster{ 0.00		}&\timeSlower{ 0.10}\\
EG\_5-12.vrep.latte		&5		&12	&\timeSlower{ 0.32	}&\timeFaster{ 0.13			}&40		&\timeFaster{ 0.15		}&\timeSlower{ 0.79}\\
EQU\_5-7a.vrep.latte		&5		&7		&\timeSlower{ 0.06		}&\timeFaster{ 0.05		}&10		&\timeFaster{ 0.01		}&\timeSlower{ 0.04}\\
EQU\_5-7b.vrep.latte		&5		&7		&\timeSlower{ 0.09		}&\timeFaster{ 0.05		}&10		&\timeSlower{ 0.01		}&\timeFaster{ 0.00}\\
HAM\_8-16.vrep.latte		&8		&16	&\timeSlower{ 1.57	}&\timeFaster{ 0.60			}&256		&\timeFaster{ 0.15		}&--\\
HC\_3-4.vrep.latte		&3		&4		&\timeSlower{ 0.03		}&\timeFaster{ 0.00		}&4		&0		&0\\
HC\_4-8.vrep.latte		&4		&8		&\timeSlower{ 0.15		}&\timeFaster{ 0.05		}&16		&\timeFaster{ 0.02		}&0.04\\
HC\_5-16.vrep.latte		&5		&16	&\timeSlower{ 0.48	}&\timeFaster{ 0.24			}&26		&\timeSlower{0.32		}&\timeFaster{ 0.25}\\
CNG\_5-6a.vrep.latte		&5		&6		&\timeSlower{ 0.04		}&\timeFaster{ 0.03		}&6		&\timeSlower{0.02		}&\timeFaster{ 0.01}\\
MJ\_32-33.vrep.latte		&32	&33	&\timeFaster{ 82.90		}&\timeSlower{83.86		}&33		&\timeSlower{1.44		}&\timeFaster{ 0.14}\\
CNG\_5-6b.vrep.latte		&5		&6		&\timeFaster{ 0.02		}&\timeSlower{ 0.03		}&6		&\timeFaster{ 0.01		}&0.02\\

\bottomrule

\end{tabular}
\end{table}

\subsection{Comparison to other software}

There are two general classes of algorithms for finding volumes and
integrals over polytopes: numerical and exact. Numerical algorithms
approximate the valuation on the polytope and involve error bounds,
whereas exact algorithms do not contain a theoretical error
term. However, exact algorithms may contain errors when they use
finite digit integers or use floating-point arithmetic. In order to
sidestep this problem, \latteInt uses NTL's arbitrary length integer and
rational arithmetic \cite{ntl-5.4} compiled with the GNU Multiple
Precision Arithmetic Library \cite{gmp-website}. The obvious downside
to exact arithmetic is speed, but this cost is necessary to obtain
exact answers. In this section, we compare our exact algorithms with
other software tools and algorithms that use numerical algorithms or non-exact
arithmetic.

{\tt Vinci} contains different algorithms for finding polytope volumes
and in fact implemented the same decompositions we used in our
software (see \cite{bueler-enge-fukuda-2000:exact-volume}). We tested
against {\tt Vinci} 1.0.5, and
Table~\ref{tabel:lawrence-vinci-volume-hot} compares {\tt LattE}'s
cone decomposition method with {\tt Vinci}'s fastest method HOT
(Hybrid Orthonormalisation Technique). We ran {\tt LattE}'s cone
decomposition method starting from the h-representation. Because the
HOT method requires both an h- and v-representation of the polytope,
we also report the time used by CDD \cite{fukuda-cdd-revisited:1996}
to convert an h-representation to a v-representation. We also break
down time spent in {\tt LattE} for finding the vertices, finding the
rays at each vertex, triangulation, and the time spent in the main
cone decomposition integration method.

%
%



\begin{table}
\centering
\caption{Time breakdown between \latteInt's cone-decomposition and \vinci's HOT method with CDD}
\label{tabel:lawrence-vinci-volume-hot}
\begin{tabular}{lrrrrrr}
	\toprule	
	 & \multicolumn{2}{c}{\vinci} & \multicolumn{4}{c}{\latte} \\
			\cmidrule(r){2-3}	\cmidrule(r){4-7}	

	Polytope	& HOT & \cdd & Vertices & Rays & Triang. & Cone. \\
	
	\midrule
	cube-9	& 0.03&0.08 & 0.02\footnotemark[1]&0.06&0.02&0.02  \\
	cube-10	& 0.11&0.18 & 0.04\footnotemark[1]&0.15&0.02&0.06 \\
	cube-14	& 141.65&7.99 & 1.24\footnotemark[1]&4.67&0.69&1.26 \\
	\hline
	rh-8-20 	& 0.13&0.89 &  0.11\footnotemark[1]&0.49&0.04&7.00 \\
	rh-8-25		& 0.43&2.63 & 0.32\footnotemark[1]&1.14&0.14&80.25 \\
	rh-10-20		& 0.96&2.21 & 0.25\footnotemark[1]&1.80&0.14&98.26 \\
	rh-10-25 	& 5.71&12.49 & 1.07\footnotemark[1]&8.80&0.44&3989.25 \\
	\hline 
	CC$_8$(9)	& 0.04&0.22 & 0.07\footnotemark[1]&0.12&0.39&0.40 \\
	CC$_8$(10)	& 0.08&0.52 & 0.16\footnotemark[1]&0.22&0.97&0.88 \\
	CC$_8$(11)	& 0.18&1.18 & 0.03\footnotemark[1]&0.42&1.84&1.76 \\
	\hline
	ccp 5		& 0.00&0.07 & 0.09\footnotemark[2]&0.00&0.10&0.09  \\
	\hline
	cross 8  & 0.00&0.39 & 0.50\footnotemark[2]&0.00&0.06&0.04 \\
	cross 9  & 0.00&1.57 & 2.15\footnotemark[2]&0.00&0.12&0.11 \\
	\hline
	rv-8-10	& 0.00&0.08 & 0.00\footnotemark[1]&0.03&0.02&0.00 \\
	rv-8-11	& 0.00&1.99 & 0.08\footnotemark[1]&0.23&0.03&0.01  \\
	rv-10-12	& 0.00&0.12 & 0.02\footnotemark[1]&0.09&0.04&0.01 \\
	rv-10-14	& 0.00&	1061.49 & 29.50\footnotemark[1]&64.96&0.10&0.07 \\
	
	\bottomrule
	\vspace{1pt} \\
	\multicolumn{7}{l}{\footnotemark[1] Computed with \fourtitwo.} \\
	\multicolumn{7}{l}{\footnotemark[2] Computed with \cdd.}	
\end{tabular}
\end{table}



It is clear that the HOT method is faster and usually accurate when applied on the \vinci database (these polytopes are available 
from \cite{bueler-enge-fukuda-2000:exact-volume}), but because of non-exact arithmetic, it can give incorrect results. In fact
we  found that {\tt Vinci}'s cone decomposition method contained a bug:  {\tt Vinci}'s cone decomposition method found the 
correct volumes for the cubes and random-hyperplane polytopes, but reported incorrect or negative volumes for most polytopes 
in the Vinci database.  We also explored how well \vinci can compute volumes of polytopes where each vertex contains small and 
large positive numbers. In Table~\ref{tabel:latte-vinci-volume-cyclic-polytopes}, we tested the accuracy of \vinci's HOT method on 
cyclic polytopes. We constructed these $d$-dimensional polytopes by taking the the convex hull of $k+d$ points 
$(t, t^2, t^3, \dots, t^d) \in \Z^d$ for $t = 5, 6, \dots, 5+k+d-1$. For very small dimensions, the HOT method does well, but gives 
incorrect or zero volumes already in dimension six. 

\begin{table}
\centering
\caption{Comparison between \latteInt and \vinci on finding the volume of cyclic polytopes}
\tabcolsep 2.0pt
\label{tabel:latte-vinci-volume-cyclic-polytopes}
\small
\begin{tabular}{*7{l}}
\toprule
    &      & \multicolumn{5}{c}{$k$} \\
\cmidrule(r){3-7}	
Dim. & Tool & 1 & 2 & 3 & 4 & 5 \\ 
 \midrule 
\multirow{2}{*}{2} & \latte & 1 & 4 & 10 & 20 & 35 \\ 
  & \vinci & 1 & 4 & 10 & 20 & 35 \\ 
 \midrule
\multirow{2}{*}{3} & \latte & 2 & 16 & 70 & 224 & 588 \\ 
  & \vinci & 1.999999999988 & 15.99999999999 & 69.99999999991 & 224.0000000006 & 587.9999999986 \\ 
 \midrule 
\multirow{2}{*}{4} & \latte & 12 & 192 & 1512 & 8064 & 33264 \\ 
  & \vinci & 11.99999993201 & 191.9999999913 & 1511.99999999 & 8063.999999892 & 33263.99999989 \\ 
 \midrule 
\multirow{2}{*}{5} & \latte & 288 & 9216 & 133056 & 1216512 & 8154432 \\ 
  & \vinci & 287.9996545868 & 9216.000252236 & 133055.9883262 & 1216511.998301 & 8154431.872519 \\ 
 \midrule
\multirow{2}{*}{6} & \latte & 34560 & 2211840 & 59304960 & 948879360 & 10600761600 \\ 
  & \vinci & 34561.951223 & 1935359.822684 & 58060819.63341 & 885910920.3761 & 10336274212.34 \\ 
 \midrule
\multirow{2}{*}{7} & \latte & 24883200 & 3185049600 & 160123392000 & 4554620928000 & 86502214656000 \\ 
  & \vinci & 25744201.0524 & 0 & 0 & 0 & 0 \\ 
 \bottomrule
\end{tabular}
\end{table}

Another comparison we made was to the paper \cite{Ong-2003:finding-exact-volume-polyhedron}, where it is claimed that 
exact volumes are computed by integration. The authors  report seven volumes for different polytopes. \latteInt's triangulation 
and cone decomposition method agrees with their calculations except in the last case. For $P_7$ the correct volume is 
$1/622080 \approx 1.607510 \times 10^{-6}$ but they calculate $1.56439\times 10^{-6}$. Presumably, because of 
non-exact arithmetic, their answer has only one digit of accuracy.

\subsection{Numerical methods}
M. Korenblit and E. Shmerling present a numerical integration algorithm in \cite{korenblit-2006:integration-over-convex-polyhedron} which is based
on a special decomposition of the integral into regions that have well-defined upper and lower limits of integration that, on an ordering of the
variables, $x_1,x_2,\dots,x_d$, $x_i$ is expressed only in terms of $x_1,\dots, x_{i-1}$.  It is known that achieving such a decomposition is
equivalent to the so-called Fourier--Motzkin elimination procedure \cite{schechter} and as such it is of exponential complexity. 
The paper \cite{korenblit-2006:integration-over-convex-polyhedron}  gives an application to finding the probability a random-coefficient polynomial has one
or two real roots in the interval $[-1,1]$. To do this, they use their software to find the volume of a polytope. They calculate $2.79167$; however, we verified that 
the correct volume is $31/12 = 2.58\bar{3}$ which gives their method one digit of accuracy.

A more interesting comparison is to \cubpack, a Fortran~90 library
which estimates the integral of a function (or vector of functions)
over a collection of $d$-dimensional hyper-rectangles and simplices
\cite{CUBPACK}. This comparison is very interesting because \cubpack
uses an adaptive grid to seek better performance and accuracy. All
integration tests with \cubpack in dimension $d$ were done with a
product of linear forms with a constant term over a random
$d$-dimensional simplex where the absolute value of any coordinate in
any vertex does not exceed 10. For example, we integrated a product of
inhomogeneous linear forms such as $(\frac{1}{5} + 2x -
\frac{37}{100}y)(2 - 5x)$ over the simplex with vertices $(10,0),
(9,9), (1,1)$.  In Table~\ref{tabel:CUBPACK-time-table}, \latte was
run 100 times to get the average running time, while \cubpack was run
1000 times due to variance.  Both the dimension and number of linear
forms multiplied to construct the integrand were varied.

\begin{table}
\centering
\caption{Average Time for \latteInt and \cubpack for integrating products of inhomogeneous linear forms over simplices.} 
\label{tabel:CUBPACK-time-table}
\tabcolsep 3pt
\footnotesize
\begin{tabular}{lcrrrrrrrrrrr}
\toprule 
 & &\multicolumn{10}{c}{Number of linear factors} \\ 
\cmidrule(r){3-12} 
$d$ & Tool & \multicolumn{1}{c}{1}  & \multicolumn{1}{c}{2}  & \multicolumn{1}{c}{3}  & \multicolumn{1}{c}{4}  & \multicolumn{1}{c}{5}  & \multicolumn{1}{c}{6}  & \multicolumn{1}{c}{7}  & \multicolumn{1}{c}{8} & \multicolumn{1}{c}{9} & \multicolumn{1}{c}{10} \\ 
\hline 
\multirow{2}{*}{2} & \latte   & \timeFaster{0.0001} & \timeFaster{0.0002} & \timeFaster{0.0005} & \timeFaster{0.0008} & \timeFaster{0.0009} & \timeFaster{0.0019} & \timeSlower{0.0038} & \timeSlower{0.0048} & \timeSlower{0.0058} & \timeFaster{0.0089} \\
                   & \cubpack & \timeSlower{0.0027} & \timeSlower{0.0014} & \timeSlower{0.0016} & \timeSlower{0.0022} & \timeSlower{0.0064} & \timeSlower{0.0052} & \timeFaster{0.0014} & \timeFaster{0.0002} & \timeFaster{0.0026} & \timeSlower{0.0213} \\
\hline 
\multirow{2}{*}{3} & \latte   & \timeFaster{0.0002} & \timeFaster{0.0005} & \timeFaster{0.0009} & \timeFaster{0.0016} & \timeFaster{0.0043} & \timeFaster{0.0073} & \timeFaster{0.0144} & \timeFaster{0.0266} & \timeFaster{0.0453} & \timeFaster{0.0748} \\
                   & \cubpack & \timeSlower{0.0134} & \timeSlower{0.0145} & \timeSlower{0.0018} & \timeSlower{0.0054} & \timeSlower{0.0234} & \timeSlower{0.0219} & \timeSlower{0.0445} & \timeSlower{0.0699} & \timeSlower{0.1170} & \timeSlower{0.2420} \\
\hline 
\multirow{2}{*}{4} & \latte   & \timeFaster{0.0003} & \timeFaster{0.0012} & \timeFaster{0.0018} & \timeSlower{0.0044} & \timeSlower{0.0121} & \timeFaster{0.0274} & \timeFaster{0.0569} & \timeFaster{0.1094} & \timeFaster{0.2247} & \timeFaster{0.4171} \\
                   & \cubpack & \timeSlower{0.0042} & \timeSlower{0.0134} & \timeSlower{0.0028} & \timeFaster{0.0019} & \timeFaster{0.0076} & \timeSlower{0.5788} & \timeSlower{4.7837} & \timeSlower{4.3778} & \timeSlower{22.3530} & \timeSlower{54.3878} \\
\hline 
\multirow{2}{*}{5} & \latte   & \timeFaster{0.0005} & \timeFaster{0.0008} & \timeTie{0.0048}    & \timeFaster{0.0108} & \timeSlower{0.0305} & \timeFaster{0.0780} & \timeFaster{0.0800} & \timeNotComputed{   0.00} & \timeNotComputed{   0.00} & \timeNotComputed{   0.00} \\
                   & \cubpack & \timeSlower{0.0013} & \timeSlower{0.0145} & \timeTie{0.0048}    & \timeSlower{0.0217} & \timeFaster{0.0027} & \timeSlower{37.0252} & \timeSlower{128.2242} & \timeNotComputed{   0.00} & \timeNotComputed{   0.00} & \timeNotComputed{   0.00} \\
\bottomrule 
\end{tabular} 
\end{table}

As shown in Table~\ref{tabel:CUBPACK-time-table}, \latteInt tends to take less time, especially when the number of forms and dimension increases.  The table does not show the high variance that \cubpack has in its run times.  For example, the 5-dimensional test case with 6 linear forms 
had a maximum running time of 2874.48 seconds, while the minimum running time was 0.05 seconds on a different random simplex.  This contrasted starkly with \latteInt, which had every test be within 0.01 (the minimum time discrepancy recognized by its timer) of every other test case.

\begin{table}
\centering
\caption{ \texttt{CUBPACK} scaling with increased relative accuracy. ``Relative Error'' is a user-specified parameter of \texttt{CUBPACK};  ``Expected Error'' is
an estimate of the absolute error, produced by \texttt{CUBPACK}'s error estimators. Finally, the ``Actual Error'' is the difference of \texttt{CUBPACK}'s result to the
exact integral computed with \latteInt.}
\label{tabel:cubpack-vs-latte}
\tabcolsep 3pt
\small
\begin{tabular}{crrrrr}
\toprule 
\multicolumn{1}{c}{Relative Error} & \multicolumn{1}{c}{Result}  & \multicolumn{1}{c}{Expected Error}  & \multicolumn{1}{c}{Actual Error}  & \multicolumn{1}{c}{\# Evaluations}  & \multicolumn{1}{c}{Time (s)}  \\ 
\midrule
$10^{-2}$ & 1260422511.762 & 9185366.414 & 94536.015 & 4467 & 0.00 \\
$10^{-3}$ & 1260507955.807 & 1173478.333 & 9091.974 & 9820 & 0.01 \\
$10^{-4}$ & 1260516650.281 & 123541.490 & 397.496 & 34411 & 0.04 \\
$10^{-5}$ & 1260517042.311 & 12588.455 & 5.466 & 104330 & 0.10 \\
$10^{-6}$ & 1260517047.653 & 1257.553 & 0.124 & 357917 & 0.31 \\
$10^{-7}$ & 1260517047.691 & 126.042 & 0.086 & 1344826 & 1.16 \\
$10^{-8}$ & 1260517047.775 & 12.601 & 0.002 & 4707078 & 4.15 \\
$10^{-9}$ & 1260517047.777 & 1.260 & $<10^{-3}$ & 16224509 & 14.09 \\
$10^{-10}$ & 1260517047.777 & 0.126 & $<10^{-3}$ & 55598639 & 48.73 \\
\bottomrule 
\end{tabular} 
\end{table}

\cubpack differs from \latteInt in that since it is based on numerical
approximations, one can ask for different levels of precision.
Table~\ref{tabel:cubpack-vs-latte} illustrates how \cubpack scales with
requested precision on a single, 4-dimensional, 10 linear form test case. It
seems that \cubpack scales linearly with the inverse of the requested
precision---10 times the precision requires about 3 times the work. 
All reported tests were done by expanding the multiplication of linear forms,
and coding a Fortran~90 function to read in the resulting polynomial and
evaluate it for specific points.  


\section{One application: Voting theory}

Computation of integrals of polynomials over polyhedral regions is
fundamental for many applications, including combinatorics,
probability and statistics. In this last section we wish to
demonstrate the power of \latteInt by attacking problems arising in
the social sciences. In the mathematical theory of voting it was
observed that the probability of events that can lead to singular
election outcomes can be modeled as the number of lattice points
inside a polytope divided by the number of lattice points of a simplex
(see \cite{ehrhart-polynomials-vothing-theory-lepelly-2007} and the
references therein).  Note that both the polytope and the simplex dilate
proportional to the number $n$ of voters.  It is very well-known from
the theory of Ehrhart functions that the counting functions are
quasipolynomials (polynomials with periodic coefficients) that depend
on $n$ \cite{barvinokzurichbook}. Thus when the quotient is evaluated
the answer is asymptotically equal to the quotient of the leading
coefficients of the two Ehrhart quasipolynomials involved.

To illustrate this, consider the following example from
\cite{ehrhart-polynomials-vothing-theory-lepelly-2007}: There are
three candidates $a, b$ and $c$, and let the preference orders of the
$n = \sum_{i=1}^6 n_i$ voters be
\begin{displaymath}
 abc \:(n_1),\: acb\: (n_2),\: bac\: (n_3),\: bca\: (n_4),\: cab\: (n_5),\: cba\: (n_6).
\end{displaymath}

Here, there are $n_1$ voters who rank candidate $a$ as first, $b$
second, and $c$ third, $n_2$ voters who rank $b$ first, $a$ second,
$c$ third, etc.  Under simple plurality voting, the candidate with the
most votes wins. But in a plurality runoff system, if no candidate
wins more than 50\% of the vote, the two candidates with the highest
vote count advance to a second voting round.
In~\cite{ehrhart-polynomials-vothing-theory-lepelly-2007}, the authors
compute the probability that the simple plurality and plurality runoff
systems give different winners. This requires setting up a system of
equations that describes the situation that $a$ wins by plurality but,
using plurality runoff $b$ obtains higher score than $c$ and a
majority of voters then prefer $b$ to $a$.

\begin{align*}
	0 & <  n_1 + n_2 - n_3 - n_4\\
	0 & <  n_3 + n_4 - n_5 - n_6, \\
	-\frac{1}{2} & <  - n_1 - n_2 - n_5,\\
	1 & =  n_1 + n_2 + n_3 + n_4 + n_5 + n_6, \\
	0 & \leq  n_i, \quad i=1,\dots, 6.	\\
\end{align*}

This is done by computing the Ehrhart quasi-polynomial of the above polyhedron and dividing by the Ehrhart quasipolynomial of
the simplex $\{\,(n_1,n_2,\dots,n_6) : n_1+n_2+\dots +n_6=1,\ n_i \geq 0\,\}$ (which is the space of all possible voting possibilities assuming
that all 6 rankings of three candidates are equally likely). All must be
multiplied by 6 because the plurality winner may be $a$, $b$ or $c$
and the second position could be $c$ not just $b$. As the authors observed, asymptotically, the leading coefficients
of these two quasipolynomials is all that matter. In the concluding remarks the authors then posed the challenge of pushing the limit of
such calculations for four-candidate elections which they observed is too big for their calculations.

However, we have observed their calculation can be further simplified and accelerated  because it is very well-known (see \cite{barvinokzurichbook}) that the  leading coefficient of the quasipolynomial is always equal to the  volume of the polytope with $n=1$, thus one can directly perform the calculation of the volume (the volume of the simplex is well-known) and do a quotient of two numbers. The key step in finding the probabilities
requires only finding the volume directly. Our algorithm corroborates that for the previous example the volume is $\frac{71}{414720}$, and when multiplied by $6 \times 120$ gives the probability these two voting systems give different winners for a large population: 12.33\%.

Using our code for exact integration we tackled the same problem for four candidates. In this case we have 24 variables associated to the
orderings

$$abcd (n_1), abdc (n_2), acbd (n_3), acdb (n_4), adbc (n_5), adcb (n_6)$$

$$bacd (n_7), badc (n_8), bcad (n_9), bcda (n_{10}), bdac (n_{11}), bdca (n_{12})$$
 
$$cabd (n_{13}), cadb (n_{14}), cbad (n_{15}), cbda (n_{16}), cdab (n_{17}), cdba (n_{18})$$

$$dabc (n_{19}), dacb (n_{20}), dbac (n_{21}), dbca (n_{22}), dcab (n_{23}), dcba (n_{24})$$

The equations and inequalities associated to the problem codify the following facts: The sum of all variables $n_i$ must be equal to the total number of
voters. We have four inequalities expressing that when $a$ is the plurality winner, $b$ obtained a score higher than~$c$, and
$c$ obtained a score higher than $d$, thus
\begin{align*}
	n_1+ n_2+ n_3+ n_4+ n_5+n_6 &> n_7+n_8+n_9+ n_{10}+n_{11}+n_{12}, \\
	n_7+n_8+n_9+n_{10}+n_{11}+n_{12} &> n_{13}+n_{14}+n_{15}+n_{16}+n_{17}+n_{18}, \\
	n_{13}+n_{14}+n_{15}+n_{16}+n_{17}+n_{18} &> n_{19}+n_{20}+n_{21}+n_{22}+n_{23}+n_{24}.
\end{align*}

These inequalities assume that the order was $a>b>c>d$ but the
answer we get should be multiplied by $4!=24$ to take into account other possible orders.
Finally, we have to express the fact that  but a majority of voters prefer $b$ over $a$ and that $a$ did not achieve more than 50 percent
of the vote ($n_1+n_2+n_3+n_4+n_5+n_6<n/2$). The volume of this polytope when $n=1$ is
\begin{displaymath}
\frac{2988379676768359}{7552997065814637134660504411827077120000}.
\end{displaymath}
The probability is then the volume times $4!$ divided by the volume of the
simplex $$\{\,(n_1,n_2,\dots,n_{24}) : \textstyle\sum_i n_i=1, \ n_i \geq
0\,\},$$ which equals $\frac{1}{23!}$. After a minute of computation using the
cone decomposition method, we obtain the probability is $12.27\%$. 

We can continue the example by considering the same problem for five candidates. The five-candidate polytope has $5!=120$ variables. However, after {\tt LRS} \cite{Avis-lrs-1999} enumerated over 12.5 million vertices, we terminated the program and decided the polytope is beyond our limits.  We close by mentioning that after the first version of this paper was made available other authors proposed new ideas to compute these values using symmetries of the problem. See \cite{achillsocialchoice}.




\section{Acknowledgements}
We are  truly grateful to our collaborators V. Baldoni, N. Berline, and M. Vergne for important discussions that led to this software.
We thank the anonymous referees for the excellent suggestions that we received from them and greatly improved the presentation of
the paper.  The senior authors, J.A. De Loera and M. K\"oppe, were partially supported by NSF grants DMS-0914107 and DMS-0914873. 
Most of the students were supported by those grants and by summer fellowships provided through the UC Davis VIGRE grant DMS-0636297.

\bibliographystyle{elsarticle-num}
\bibliography{biblio}

\begin{thebibliography}{10}
\expandafter\ifx\csname url\endcsname\relax
  \def\url#1{\texttt{#1}}\fi
\expandafter\ifx\csname urlprefix\endcsname\relax\def\urlprefix{URL }\fi
\expandafter\ifx\csname href\endcsname\relax
  \def\href#1#2{#2} \def\path#1{#1}\fi

\bibitem{BLi:1993}
B.~Li, The moment calculation of polyhedra, Pattern Recognition 26~(8) (1993)
  1229--1233.

\bibitem{Mirtich:1996}
B.~Mirtich, Fast and accurate computation of polyhedral mass properties,
  Journal of Graphics Tools. 1~(2) (1996) 31--50.

\bibitem{Ong-2003:finding-exact-volume-polyhedron}
H.~Ong, H.~Huang, W.~Huin, Finding the exact volume of a polyhedron, Advances
  in Engineering Software 34~(6) (2003) 351 -- 356.
\newblock \href {http://dx.doi.org/DOI: 10.1016/S0965-9978(03)00030-9}
  {\path{doi:DOI: 10.1016/S0965-9978(03)00030-9}}.

\bibitem{Sturmfels1996}
B.~Sturmfels, Gr\"obner bases and convex polytopes, Vol.~8 of University
  Lecture Series, American Mathematical Society, 1996.

\bibitem{Berline-Vergne-2007}
N.~Berline, M.~Vergne, Local {E}uler--{M}aclaurin formula for polytopes, Moscow
  Math. J. 7 (2007) 355--386.

\bibitem{Bronstein}
M.~Bronstein, Symbolic Integration {I} -- Transcendental Functions, Vol.~1,
  Springer, Heidelberg, 2005.

\bibitem{dyerfrieze88}
M.~E. Dyer, A.~M. Frieze, On the complexity of computing the volume of a
  polyhedron, SIAM J. Comput. 17~(5) (1988) 967--974.

\bibitem{brightwellwinkler91}
G.~Brightwell, P.~Winkler, Counting linear extensions, Order 8~(3) (1991)
  225--242.

\bibitem{khachiyan93}
L.~Khachiyan, Complexity of polytope volume computation, in: New trends in
  discrete and computational geometry, Vol.~10 of Algorithms Combin., Springer,
  Berlin, 1993, pp. 91--101.

\bibitem{lawrence91}
J.~Lawrence, Polytope volume computation, Math. Comp. 57~(195) (1991) 259--271.

\bibitem{gritzmann-complexity-volume:1994}
P.~Gritzmann, V.~Klee, On the complexity of some basic problems in
  computational convexity: {II}. {V}olume and mixed volumes, Universit{\"a}t
  Trier, Mathematik/Informatik, Forschungsbericht 94-07.

\bibitem{elekes86}
G.~Elekes, A geometric inequality and the complexity of computing volume,
  Discrete Comput. Geom. 1~(4) (1986) 289--292.

\bibitem{rademacher}
L.~Rademacher, Approximating the centroid is hard, in: Proceedings of 23th
  annual ACM Symposium of Computational Geometry, Gyeongju, South Korea, June
  6-8, 2007, 2007, pp. 302--305.

\bibitem{latteintegrale}
J.~De~Loera, B.~Dutra, M.~K{\"o}ppe, S.~Moreinis, G.~Pinto, J.~Wu, A users
  guide for latte integrale v1.5, Available from URL
  {\url{http://www.math.ucdavis.edu/~latte/}} (2011).

\bibitem{baldoni-berline-deloera-koeppe-vergne:integration}
V.~Baldoni, N.~Berline, J.~A. De~Loera, M.~K{\"o}ppe, M.~Vergne, How to
  integrate a polynomial over a simplex, Math. Comp. 80~(273) (2011) 297--325.
\newblock \href {http://dx.doi.org/10.1090/S0025-5718-2010-02378-6}
  {\path{doi:10.1090/S0025-5718-2010-02378-6}}.

\bibitem{bueler-enge-fukuda-2000:exact-volume}
B.~B{\"u}eler, A.~Enge, K.~Fukuda, Exact volume computation for polytopes: A
  practical study, in: G.~Kalai, G.~M. Ziegler (Eds.), Polytopes --
  Combinatorics and Computation, Vol.~29 of DMV-Seminars, Birkh{\"a}user
  Verlag, Basel, 2000.

\bibitem{Barvinok-1991}
A.~I. Barvinok, Computation of exponential integrals, Zap. Nauchn. Sem.
  Leningrad. Otdel. Mat. Inst. Steklov. (LOMI) Teor. Slozhn. Vychisl. 5 (1991)
  149--162, 175--176, translation in J. Math. Sci. 70 (1994), no. 4,
  1934--1943.

\bibitem{Barvinok-1992}
A.~I. Barvinok, Partition functions in optimization and computational problems,
  Algebra i Analiz 4 (1992) 3--53, translation in St. Petersburg Math. J. 4
  (1993), no. 1, pp. 1--49.

\bibitem{Brion88}
M.~Brion, Points entiers dans les poly{\`e}dres convexes, Ann. Sci. {\'E}cole
  Norm. Sup. 21~(4) (1988) 653--663.

\bibitem{lasserre-algorithm-volume-polyhedron:1983}
J.~B. Lasserre, An analytical expression and an algorithm for the volume of a
  convex polyhedron in $\mathbf{R}^n$, Journal of Optimization Theory and
  Applications 39 (1983) 363--377, 10.1007/BF00934543.
\newblock \href {http://dx.doi.org/10.1007/BF00934543}
  {\path{doi:10.1007/BF00934543}}.

\bibitem{brion-lawrence-varchenko-rational-generating-functions-2005}
M.~Beck, C.~Haase, F.~Sottile, Formulas of {B}rion, {L}awrence, and {V}archenko
  on rational generating functions for cones, The Mathematical Intelligencer 31
  (2009) 9--17, 10.1007/s00283-008-9013-y.
\newblock \href {http://dx.doi.org/10.1007/s00283-008-9013-y}
  {\path{doi:10.1007/s00283-008-9013-y}}.

\bibitem{barvinokzurichbook}
A.~I. Barvinok, Integer Points in Polyhedra, Z\"urich Lectures in Advanced
  Mathematics, European Mathematical Society (EMS), Z\"urich, Switzerland,
  2008.

\bibitem{BarviPom}
A.~I. Barvinok, J.~E. Pommersheim, An algorithmic theory of lattice points in
  polyhedra, in: L.~J. Billera, A.~Bj\"orner, C.~Greene, R.~E. Simion, R.~P.
  Stanley (Eds.), New Perspectives in Algebraic Combinatorics, Vol.~38 of Math.
  Sci. Res. Inst. Publ., Cambridge Univ. Press, Cambridge, 1999, pp. 91--147.

\bibitem{beckrobins}
M.~Beck, S.~Robins, Computing the continuous discretely: integer-point
  enumeration in polyhedra, Undergraduate Texts in Mathematics, Springer, 2007.

\bibitem{Avis-lrs-1999}
D.~Avis, lrs: A revised implementation of the reverse search vertex enumeration
  algorithm, \url{http://cgm.cs.mcgill.ca/~avis/C/lrs.html} (1999).

\bibitem{fukuda-cdd-revisited:1996}
K.~Fukuda, A.~Prodon, Double description method revisited, in: M.~Deza,
  R.~Euler, I.~Manoussakis (Eds.), Combinatorics and Computer Science, Vol.
  1120 of Lecture Notes in Computer Science, Springer Berlin / Heidelberg,
  1996, pp. 91--111.
\newblock \href {http://dx.doi.org/10.1007/3-540-61576-8_77}
  {\path{doi:10.1007/3-540-61576-8_77}}.

\bibitem{alexanderhirschowitz}
J.~Alexander, A.~Hirschowitz, Polynomial interpolation in several variables, J.
  Algebraic Geom. 4 (1995) 201--222.

\bibitem{brambillaottaviani}
M.~C. Brambilla, G.~Ottaviani, On the {A}lexander--{H}irschowitz theorem,
  Journal of Pure and Applied Algebra 212~(5) (2008) 1229--1251, available at
  arXiv:math.AG/0701409v2.

\bibitem{Carlinietal}
E.~Carlini, M.~Catalisano, A.~Geramita, The solution to {W}aring's problem for
  monomials, eprint arXiv:{\penalty0}math.AC/1110.0745v1 (2011).

\bibitem{Henrici}
P.~Henrici, Applied and computational complex analysis. {V}ol. 1, Wiley
  Classics Library, John Wiley \& Sons Inc., New York, 1988, power
  series---integration---conformal mapping---location of zeros, Reprint of the
  1974 original, A Wiley-Interscience Publication.

\bibitem{barvinok-woods-2003}
A.~I. Barvinok, K.~Woods, Short rational generating functions for lattice point
  problems, Journal of the AMS 16~(4) (2003) 957--979.

\bibitem{latte1}
J.~A. De~Loera, R.~Hemmecke, J.~Tauzer, R.~Yoshida, Effective lattice point
  counting in rational convex polytopes, Journal of Symbolic Computation 38~(4)
  (2004) 1273--1302.

\bibitem{koeppe:irrational-barvinok}
M.~K{\"o}ppe, A primal {B}arvinok algorithm based on irrational decompositions,
  SIAM Journal on Discrete Mathematics 21~(1) (2007) 220--236.
\newblock \href {http://dx.doi.org/10.1137/060664768}
  {\path{doi:10.1137/060664768}}.

\bibitem{trip-burst-tries-gastineaau-2006}
M.~Gastineau, J.~Laskar, Development of {TRIP}: Fast sparse multivariate
  polynomial multiplication using burst tries, in: V.~Alexandrov, G.~van
  Albada, P.~Sloot, J.~Dongarra (Eds.), Computational Science -- ICCS 2006,
  Vol. 3992 of Lecture Notes in Computer Science, Springer Berlin / Heidelberg,
  2006, pp. 446--453.
\newblock \href {http://dx.doi.org/10.1007/11758525_60}
  {\path{doi:10.1007/11758525_60}}.

\bibitem{cddlib-094a}
K.~Fukuda, {\texttt{cddlib}}, version 094f, Available from URL
  \url{http://www.ifor.math.ethz.ch/~fukuda/cdd_home/} (2008).

\bibitem{4ti2}
4ti2 team, 4ti2---a software package for algebraic, geometric and combinatorial
  problems on linear spaces, {A}vailable at \url{www.4ti2.de}.

\bibitem{Guillaumeetal}
G.~Hanrot, P.~Zimmermann, {A long note on Mulders' short product}, Rapport de
  recherche RR-4654, INRIA, available at
  \url{http://hal.inria.fr/inria-00071931} (2002).

\bibitem{CUBPACK}
R.~Cools, A.~Haegemans, Algorithm 824: {CUBPACK}: a package for automatic
  cubature; framework description., ACM Trans. Math. Software 29~(3) (2003)
  287--296.

\bibitem{korenblit-2006:integration-over-convex-polyhedron}
M.~Korenblit, E.~Shmerling, Algorithm and software for integration over a
  convex polyhedron, in: A.~Iglesias, N.~Takayama (Eds.), Mathematical Software
  -- ICMS 2006, Vol. 4151 of Lecture Notes in Computer Science, Springer
  Berlin, Heidelberg, 2006, pp. 273--283.
\newblock \href {http://dx.doi.org/10.1007/11832225_28}
  {\path{doi:10.1007/11832225_28}}.

\bibitem{polymake-software}
E.~Gawrilow, M.~Joswig, polymake: a framework for analyzing convex polytopes,
  in: G.~Kalai, G.~M. Ziegler (Eds.), Polytopes --- Combinatorics and
  Computation, Birkh\"auser, 2000, pp. 43--74.

\bibitem{ziegler-db-website}
G.~Ziegler, Database of 0--1 polytopes,
  \url{http://www.math.tu-berlin.de/polymake/examples/}, accessed: 07/19/2011.

\bibitem{extratables}
J.~De~Loera, B.~Dutra, M.~K{\"o}ppe, S.~Moreinis, G.~Pinto, J.~Wu, Software for
  exact integration of polynomials over polyhedra: online supplement, Available
  from URL
  {\url{http://www.math.ucdavis.edu/~latte/theory/SoftwareExactIntegrationPoly%
nomialsPolyhedraOnlineSupplement.pdf}} (2012).

\bibitem{ntl-5.4}
V.~Shoup, {NTL}, a library for doing number theory, Available from URL
  \url{http://www.shoup.net/ntl/} (2005).

\bibitem{gmp-website}
The {GNU} multiple precision arithmetic library, \url{http://gmplib.org/}.

\bibitem{schechter}
M.~Schechter, Integration over a polyhedron: An application of the
  {F}ourier--{M}otzkin elimination method, American Mathematical Monthly
  105~(3) (1998) 246--251.

\bibitem{ehrhart-polynomials-vothing-theory-lepelly-2007}
D.~Lepelley, A.~Louichi, H.~Smaoui, On {E}hrhart polynomials and probability
  calculations in voting theory, Social Choice and Welfare 30 (2008) 363--383.
\newblock \href {http://dx.doi.org/10.1007/s00355-007-0236-1}
  {\path{doi:10.1007/s00355-007-0236-1}}.

\bibitem{achillsocialchoice}
A.~Sch\"urmann, Exploiting polyhedral symmetries in social choice, Available
  from URL \url{http://front.math.ucdavis.edu/1109.1545} (2012).

\end{thebibliography}
\end{document}